\documentclass[preprint,12pt]{elsarticle}

\usepackage{amsmath}
\usepackage{array,overpic}
\usepackage{epsf,latexsym,amsfonts,amsbsy}
\usepackage{amsmath, amssymb, amsthm,amscd,amsxtra}
\usepackage{graphicx,epstopdf}
\usepackage{color}

\def\be {\begin{equation}}
\def\ee {\end{equation}}
\def\ba {\begin{eqnarray}}
\def\ea {\end{eqnarray}}

\newcommand{\tss}{\textstyle}
\newtheorem{theorem}{Theorem}
\newtheorem{lemma}[theorem]{Lemma}
\newdefinition{rmk}{Remark}

\usepackage{amssymb}
\journal{J. Differential Equation}

\begin{document}

\begin{frontmatter}



\title{Bifurcation of ten small-amplitude limit cycles \\
by perturbing a quadratic Hamiltonian system \\
with cubic polynomials}


\author{Yun Tian, Pei Yu\footnote{Corresponding author.\\
\hspace*{0.25in}E-mail addresses: ytian56@uwo.ca (Y. Tian), pyu@uwo.ca (P. Yu)}}

\address{Department of Applied Mathematics, Western University  \\
London, Ontario, Canada N6A 5B7}

\begin{abstract}
This paper contains two parts. In the first part,
we shall study the Abelian integrals for \.Zo{\l}\c{a}dek's
example~\cite{Zoladek1995}, in which it is claimed the existence
of $11$ small-amplitude limit cycles around
a singular point in a particular cubic vector filed.
We will show that the basis chosen in the proof of~\cite{Zoladek1995}
are not independent, which leads to failure in drawing the conclusion of
the existence of $11$ limit cycles in this example.
In the second part, we present a good combination of Melnikov function
method and focus value (or normal form) computation method to study
bifurcation of limit cycles.
An example by perturbing a quadratic
Hamiltonian system with cubic polynomials is presented to demonstrate the
advantages of both methods, and $10$ small-amplitude limit cycles
bifurcating from a center are obtained by using up to 5th-order Melnikov
functions.

\end{abstract}


\end{frontmatter}

% \linenumbers

\section{Introduction}
The well-known Hilbert's 16th problem~\cite{Hilbert1900} has been studied
for more than one century, and the research on this problem is still very
active today. To be more specific, consider the following planar system:
\begin{eqnarray}
\dot{x} = P_n(x,y), \qquad \dot{y} = Q_n(x,y),
\label{Eq1}
\end{eqnarray}
where $\, P_n(x,y) \,$ and $\, Q_n(x,y) \,$ represent $n$th-degree
polynomials in $\, x\,$ and $\, y$. The second part of Hilbert's 16th
problem is to find the upper bound, called Hilbert number $\, H(n) $,
on the number of limit cycles that system (\ref{Eq1}) can have.

The progress in the solution of the problem is very slow. Even the simplest
case $n=2$ has not been completely solved, though in the early 1990's,
Ilyashenko~\cite{Ilyashenko1991} and \'{E}calle~\cite{Ecalle1992}
independently proved that the number of limit cycles is finite
for any given planar polynomial vector field. For general quadratic
polynomial systems, the best result is $ H(2) \ge 4$, obtained more than
$30$ years ago~\cite{Shi1980,ChenWang1979}. Recently, this result was also
obtained for near-integrable quadratic systems~\cite{YH2012}.
However, whether $\, H(2)=4 \,$ is still open. For cubic polynomial systems,
many results have been obtained on the lower bound of the Hilbert number.
So far, the best result for cubic systems is
$\, H(3) \ge 13$~\cite{LiLiu2010,LLY09}. Note that the $13$ limit cycles
obtained in~\cite{LiLiu2010,LLY09} are distributed around several
singular points.
%This number is believed to be below the maximal number
%which can be obtained for generic cubic systems.
A comprehensive review
on the study of Hilbert's 16th problem can be found in a survey
article~\cite{Li2003}.

%\pagestyle{myheadings}
%\markright{\centerline{\footnotesize {\it Y. Tian, P. Yu/
%J. Differential Equations {\rm xxx} (201{\rm x}) {\rm xxx-xxx}}}}
\renewcommand{\leftmark}{Chapter\ \thechapter. \ Ten limit cycles in a cubic near-Hamiltonian system}

In order to help understand and attack Hilbert's 16th problem, the so
called weakened Hilbert's 16th problem was posed by Arnold~\cite{Arnold83}.
The problem is to ask for the maximal number of isolated zeros of the
Abelian integral or Melnikov function:
\begin{equation}
M(h,\delta) = \displaystyle\oint_{H(x,y)=h} Q(x,y) \, dx - P(x,y) \, dy,
\label{Eq2}
\end{equation}
where $\, H(x,y), \, P(x,y) \,$ and $\, Q(x,y) \,$ are all real polynomials
in $\,x \,$ and $\, y\,$, and the level curves $H(x,y)=h$ represent at least
a family of closed orbits for $h\in(h_1,h_2)$, and $\delta$ denotes the
parameters (or coefficients) involved in $P$ and $Q$. The weakened Hilbert's
16th problem itself is a very important and interesting problem, closely
related to the study of limit cycles in the following near-Hamiltonian
system~\cite{Han2006}:
\begin{equation}
\dot{x} = H_y(x,y) + \varepsilon \, P(x,\,y), \qquad
\dot{y} =- \,  H_x(x,y) + \varepsilon \, Q(x,\,y),
\label{Eq3}
\end{equation}
where $0 < \varepsilon \ll 1$. Studying the
bifurcation of limit cycles for such a system can be now transformed to
investigating the zeros of the Melnikov function $M(h,\delta)$,
given in (\ref{Eq2}).

When we focus on the maximum number of small-amplitude limit cycles,
denoted by $M(n)$, bifurcating from an elementary center or an elementary focus,
one of the best-known results is $M(2)=3$, which was proved by Bautin
in 1952~\cite{Bautin1952}. For $n=3$, several results have been obtained
(e.g. see \cite{Zoladek1995,YuHan2011,James1991}). Among them, in 1995
\.Zo{\l}\c{a}dek~\cite{Zoladek1995} first constructed a rational Darboux
integral to study existence of $11$ small-amplitude limit cycles in
cubic vector fields.
This pioneer work later motivated many researches in this area to study
bifurcation of limit cycles. The rational
Darboux integral in~\cite{Zoladek1995} is given by
\begin{equation}\label{Eq4}
H_0=\frac{f_1^5}{f_2^4}=\frac{(x^4+4x^2+4y)^5}{(x^5+5x^3+5xy+5x/2+a)^4},
\end{equation}
which generates the dynamical system in the form of
\begin{equation}\label{Eq5}
\begin{split}
\dot x&= x^3+xy+5x/2+a,\\
\dot y&= -ax^3+6x^2y-3x^2+4y^2+2y-2ax,
\end{split}
\end{equation}
with the integrating factor $M=20f_1^4f_2^{-5}$.

It can be shown that for $a<-2^{5/4}$, system (\ref{Eq5})
has a center $C_0=(-a/2,-a^2/4-1/2)$
and five (real or complex) critical points $(r,-r^2-5/2-a/r)$, where $\,r\,$
satisfies the polynomial equation $r^5-10r-4a=0$. In addition,
system (\ref{Eq5}) has a saddle point and a non-elementary critical
point at infinity. Let $h_0=H_0(C_0)=-2/a$. Around $C_0$, there exists a
family of periodic orbits given by $\{\gamma_h: H_0(x,y)=h$,
$0<h-h_0\ll1\}$. $\gamma_{h}$ approaches $C_0$ as $h\rightarrow h_0^+$.

Recently, Yu \& Han~\cite{YuHan2011} applied a different method to study
system (\ref{Eq5}) and only got $9$ small-amplitude limit cycles around
the center $C_0$.
This difference motivated us to reconsider system (\ref{Eq5}) carefully
and finally find some not-easy-to-find mistakes in the proof
of~\cite{Zoladek1995}, which lead to failure in drawing the conclusion
of existence of $11$ limit cycles.
In the next section, we shall present a detailed analysis on
the Abelian integrals of system (\ref{Eq5})
and point out where the mistakes were made in~\cite{Zoladek1995}.
However, we must emphasize that although some flaws were found
in~\cite{Zoladek1995}, the idea and methodology
presented in this paper are still very valuable and useful.
In fact, our work was motivated by this paper as well as that~\cite{Iliev2000}.

In the second part, we will present a good combination of
higher-order Melnikov function method and focus value computation method
to study the number of small-amplitude limit cycles. As a matter of fact,
in proving existence of small-amplitude limit cycles,
higher-order Melnikov functions are equivalent to higher $\varepsilon$-order
focus values. These two methods have their advantages and disadvantages:
Melnikov function method can be used to show the vanishing of a particular
order Melinikov function, which is equivalent to show
the vanishing of all this particular $\varepsilon$-order focus values, and this
cannot be proved by using focus value computation since one cannot verify
an infinite number of focus values. On the other hand,
for higher-order Melnikov functions, it becomes extremely difficult
to prove the independence of a choosing set of basis,
while this is straightforward by using focus value
computation. In particular, we present an example to
demonstrate this good combination method and
obtain $10$ small-amplitude limit
cycles by perturbing a quadratic Hamiltonian system with 3rd-degree
polynomial functions.

In general, a perturbed quadratic Hamiltonian
system can be described by
\begin{equation}\label{Eq6}
\begin{split}
& \dot x= y+a_1xy+a_2y^2+\varepsilon P(x, y, \varepsilon),\\
& \dot y= -x+x^2-\frac{1}{2}a_1y^2+\varepsilon Q(x, y, \varepsilon),
\end{split}
\end{equation}
where $P$ and $Q$ are $n$th-degree polynomials in $x$ and $y$ with
coefficients depending analytically on the small parameter $\varepsilon$.
When $\varepsilon=0$, system (\ref{Eq6}) has a cubic Hamiltonian,
\begin{equation}\label{Eq7}
H(x, y)=\frac{1}{2}(x^2+y^2)-\frac{1}{3}x^3+\frac{1}{2}a_1xy^2
+\frac{1}{3}a_2y^3,
\end{equation}
and its parameters $a_1$ and $a_2$ take values from the set,
\begin{equation*}
\Omega=\left\{(a_1,a_2)\in \mathbb{R}^2:
-1\le a_1\le 2,\ 0\le a_2\le (1-\frac{a_1}{2}) \sqrt{1+a_1}
\right\}.
\end{equation*}
The Hamiltonian given in (\ref{Eq7}) is actually the so-called normal
form~\cite{HorozovIliev1994} for all quadratic Hamiltonian systems which
have a center at the origin. There exists a family of closed ovals around
the origin given by $\{\Gamma_h: H(x,y)=h$, $h\in(0,\frac{1}{6})\}$.

For any $h\in(0,\frac{1}{6})$ the displacement function $d(h,\epsilon)$
of system (\ref{Eq6}) has a representation
\begin{equation}\label{Eq8}
d(h,\varepsilon)=\varepsilon M_1(h)+\varepsilon^2M_2(h)
+\varepsilon^3M_3(h)+\cdots,
\end{equation}
where $M_i(h)$ is called the $i$th-order Melnikov function, particularly
the higher-order Melnikov functions if $i\ge 2$. Then, we may determine
the number of the limit cycles of system (\ref{Eq6}) emerging from the
closed ovals $\{\Gamma_h\}$ by studying the zeros of the first
non-vanishing Melnikov function $M_i(h)$ in $h\in(0,\frac{1}{6})$.

Suppose $M_1(h)\not\equiv0$ in (\ref{Eq8}). Denote $Z(n)$ the sharp upper
bound of the number of zeros of $M_1(h)$ for system (\ref{Eq6}), where
$n=\max\{\deg(P),\deg(Q)\}$. Gavrilov~\cite{Gavrilov2001}
proved $Z(2)=2$ for the Hamiltonian $H$ with four distinct critical
values (in a complex domain). Horozov \& Iliev~\cite{HorozovIliev1998}
gave a linear estimate $Z(n)\le 5 (n+3)$. Also, some sharp upper bounds
were obtained for certain particular cubic Hamiltonians, for example:
$n\!-\!1$ for the Bogdanov-Takens Hamiltonian,
$H=\frac{1}{2}(x^2+y^2)-\frac{1}{3}x^3$
(see~\cite{Petrov1986}), and $[\frac{2}{3}(n-1)]$ 
for the Hamiltonian triangle,
$H=\frac{1}{2}(x^2+y^2)-\frac{1}{3}x^3+xy^2$ (see~\cite{Gavrilov1998}).

Moreover, for the Bogdanov-Takens Hamiltonian, there are some results on
the upper bound of the number of zeros of the first nonvanishing higher-order
Melnikov function $M_k(h)$. Li \& Zhang~\cite{LiZhang1995} got a sharp
upper bound for $k=2$: $2n-2$ when $n$ is even and $2n-3$ when $n$ is odd.
Iliev~\cite{Iliev2000} obtained a sharp upper bound $3n-4$ for $k=3$,
by applying the Fran\c{c}oise's procedure~\cite{Francoise1996} for computing
higher-order Melnikov functions. The higher-order Melnikov functions can be
also easily expressed via iterated integrals, which will be seen in the
next section.

In this paper, we study the number of small-amplitude limit cycles
in (\ref{Eq6}) bifurcating from the origin, using higher-order Melnikov
functions. Hereafter we suppose $P$ and $Q$ are cubic polynomials
in the form of
\begin{equation}\label{Eq9}
\begin{split}
&P(x,y,\varepsilon)=\sum_{m=1}^\infty\varepsilon^{m-1}P_m(x,y)\quad
\mbox{with}\ \ P_m(x,y)=\sum_{i+j=1}^{3}a_{ijm}x^iy^j,\\
&Q(x,y,\varepsilon)=\sum_{m=1}^\infty\varepsilon^{m-1}Q_m(x,y)\quad
\mbox{with}\ \ Q_m(x,y)=\sum_{i+j=1}^{3}b_{ijm}x^iy^j.
\end{split}
\end{equation}

Our main result is given below, and its proof will be given in
Section \ref{proof_Thm1}.

\noindent
\begin{theorem}\label{TheoB1}
Let the polynomials $P$ and $Q$ in (\ref{Eq6}) be
given by (\ref{Eq9}). Then for any $1\le k\le 5$,
there exist real values for $(a_1,a_2)\in\Omega$ such that
system $(\ref{Eq6})$ can have $\lfloor\frac{4k}{3}\rfloor\!+\!4$
small-amplitude limit cycles around the origin under proper
cubic perturbations, when $M_k(h)$ is the first
non-vanishing Melnikov function in $(\ref{Eq8})$.
\end{theorem}

%\noindent{\bf Remark 1.}
\begin{rmk}
It follows from Theorem \ref{TheoB1} that $10$ small-amplitude limit
cycles exist in the vicinity of the origin of system (\ref{Eq6}) when $k=5$.
\end{rmk}

The rest of the paper is organized as follow. In the next section,
we consider the example given in~\cite{Zoladek1995},
and show that the bases chosen in the proof
are not independent, leading to failure in drawing the conclusion of
the existence of $11$ limit cycles.
In Section \ref{Hamiltonian}, we present some results for
polynomial one-forms with respect to the Hamiltonian (\ref{Eq7}), which are
needed for the proof of Theorem 1 in Section \ref{proof_Thm1}.
Then, in Section \ref{proof_Thm1} by choosing special
forms for the polynomials $P$ and $Q$ without loss of generality, we prove
Theorem \ref{TheoB1}. Finally, conclusion is drawn in Section \ref{concl}.

\section{Abelian integrals of system (\ref{Eq5})}\label{exam}

In this section, we consider system (\ref{Eq5})
and briefly describe the methodology used in \cite{Zoladek1995}.
Suppose the perturbed system of (\ref{Eq5}) is
described by
\begin{equation}\label{Eq10}
\begin{split}
\dot x &= M^{-1}H_{0y} + \varepsilon p(x,y,\varepsilon),\\
\dot y &= -M^{-1}H_{0x} + \varepsilon q(x,y,\varepsilon),
\end{split}
\end{equation}
where $p(x,y,\varepsilon)$ and $q(x,y,\varepsilon)$ are polynomials in $x$
and $y$ with coefficients depending analytically on the small parameter
$\varepsilon$, and $\max(\deg(p),\deg(q))\le 3$.
Note that the non-perturbed system (\ref{Eq10})
(i.e. $\varepsilon=0$) has a center at $C_0$.

Let $S$ be a section transversal to the closed orbit
$\{\gamma_h: H_0(x,y)=h,\ 0<h-h_0\ll1\}$, with $h$
as a parameter, we define the Poincar\'e map
$\mathcal{P}(h,\varepsilon)$ of system (\ref{Eq10}), and thus the
corresponding displacement function,
$d(h,\varepsilon)=\mathcal{P}(h,\varepsilon)-h$, has the form
\begin{equation}\label{Eq11}
d(h,\varepsilon)=\varepsilon \int_{L(h,\varepsilon)}M(qdx-pdy)
=\varepsilon M_1(h) +\varepsilon^2 M_2(h) + O(\varepsilon^3),
\end{equation}
where $L(h,\varepsilon)$ is a trajectory of the perturbed system (\ref{Eq10}).
We can use the first non-vanishing Melnikov function $M_k(h)$ in (\ref{Eq11})
to investigate the number of the limit cycles around the center $C_0$.
Generally, the zeros of $M_k(h)$ correspond to the limit
cycles of system (\ref{Eq10}) around $C_0$.

Let $\varpi=qdx-pdy$, $\deg(\varpi)=\max(\deg(p),\deg(q))$. Then, the
first-order Melnikov function $M_1(h)$ can be expressed in the form of
$$
M_1(h)=\oint_{\gamma_h}M\varpi|_{\varepsilon=0}
=h\oint_{\gamma_h} \left. \frac{\varpi}{f_1f_2} \right|_{\varepsilon=0}.
$$
When $M_1(h)\equiv0$, we may use an iterated integral to express the
second-order Melnikov function $M_2(h)$. The first integral of system
(\ref{Eq10}) can be approximated as $H_{\varepsilon}=H_0-\varepsilon H_1$,
where the function $H_1$ is defined by
$H_1(B)=\int_A^BM\varpi|_{\varepsilon=0}$, evaluated along the orbit $\gamma_h$,
with $A=\gamma_h\bigcap S$ and $B\in\gamma_h$. Thus, for system (\ref{Eq10})
we have the second-order Melnikov function, given by
\begin{equation}\label{Eq12}
M_2(h)= \left. \frac{d}{d\varepsilon}\Big(\int_{H_\varepsilon=h}
M\varpi\Big) \right|_{\varepsilon=0}.
\end{equation}
Suppose that the polynomials $p$ and $q$ are expanded as
$$
\begin{array}{ll}
p(x,y,\varepsilon)=p_1(x,y)+\varepsilon p_2(x,y)+O(\varepsilon^2),\\[0.5ex]
q(x,y,\varepsilon)=q_1(x,y)+\varepsilon q_2(x,y)+O(\varepsilon^2).
\end{array}
$$
Further, let $\varpi_i=q_idx-p_idy$, $i=1,2$. Then (\ref{Eq12}) can be
rewritten as
\begin{equation}\label{Eq13}
\begin{split}
M_2(h)&= \left. \frac{d}{d\varepsilon}\Big(\int_{H_\varepsilon =h}
M\varpi_1\Big) \right|_{\varepsilon=0} +\oint_{\gamma_h}M\varpi_2
      =\oint_{\gamma_h}\frac{d(M\varpi_1)}{dH_0}H_1
+h\oint_{\gamma_h}\frac{\varpi_2}{f_1f_2},
\end{split}
\end{equation}
where $\frac{d(M\varpi_1)}{dH_0}$ is the Gelfand-Leray form
(see \cite{Arnold1988}).

In~\cite{Zoladek1995}, the second-order Melnikov function $M_2(h)$
was used to study the small-amplitude limit cycles
of system (\ref{Eq10}) bifurcating from the center $C_0$.
More precisely, twelve Abelian integrals
$I_{\omega_i}(h)=\oint_{\gamma_h}\omega_i/(f_1f_2)$, $i=1,\ldots,12$,
were chosen for (\ref{Eq13}),
where one-forms $\omega_i$ are given as follows:
$$
\begin{array}{ll}
\omega_k= x^{k-1} dx,\, k=1,2,3,4, \quad
\omega_5=(18x^2+18y+5)dx, \quad  \omega_6=xydx, \\[1.0ex]
\omega_7=x^2ydx, \quad \omega_8=xy^2dx, \quad
\omega_9=y^3dx, \quad \omega_{10}=xy^2dy, \quad \omega_{11}=y^3dy, \\ [1.0ex]
\omega_{12}=y^2(5-3x^2)dx+xy(x^2+1)dy.
\end{array}
$$
Then, by showing the independence of the integrals
$I_{\omega_i}(h)$, $1\le i\le 12$,
it is claimed in~\cite{Zoladek1995}
that $11$ small-amplitude limit cycles can bifurcate
from the center $C_0$ after suitable cubic perturbations.

Later, system (\ref{Eq10}) was re-investigated by using the method of focus
values computation~\cite{YuHan2011}. Based on the computation of
$\varepsilon$-order and $\varepsilon^2$-order focus values, the authors
of~\cite{YuHan2011} showed that system (\ref{Eq10}) could have $9$
small-amplitude limit cycles bifurcating from the center $C_0$. This obvious
difference motivated us to study system (\ref{Eq5}), and finally
to find that any vector of the linear space of integrals
$I_{\omega}(h)=\oint_{\gamma_h}\omega/(f_1f_2)$, $\deg(\omega)\le3$,
can be expressed as a linear combination of the ten
integrals $I_{\omega_i}(h)$, $1\le i\le 11$, $i\neq 4$.
In the following, we show the details.

Nine one-forms $\eta_j$, $1\le j\le 9$, were obtained in \cite{Zoladek1995}
satisfying $I_{\eta_j}(h)=0$, where
\begin{equation*}
\begin{split}
\eta_1&= (x^3+2x)dx+dy,\\
\eta_2&= (-3ax^2+12xy-6x-2a)dx-(3x^2+y+ \tss\frac{5}{2})dy,\\
\eta_3&= (6x^2+8y+2)dx-xdy,\\
\eta_4&= (-3ax^3+12x^2y-6x^2-2ax)dx-(2x^3-a)dy,\\
\eta_5&= (ax^3+3x^2+4y^2+2ax)dx-xydy,\\
\eta_6&= (-ax^3+6x^2y-3x^2+4y^2+2y-2ax)dx-(x^3+xy+ \tss\frac{5}{2} x+a)dy, \\
\eta_7&= (3ax^2y \!-\! 12xy^2 \!+\! 6xy \!+\! 2ay)dx
-(3x^2y \!-\! ax^3 \!-\! 3x^2 \!+\! 3y^2 \!-\! \tss\frac{1}{2} y
\!-\! 2ax)dy,\\
\eta_8&= (-5x^3-7xy+ \tss\frac{1}{2} x+a)dx+x^2dy,\\
\eta_9&= (\tss\frac{21}{2} xy-7xy^2+ay)dx+(2x^2y-\tss\frac{3}{2} x^2+ax+y)dy.
\end{split}
\end{equation*}
We find another one-form $\eta_{10}$, given by
$$
\eta_{10}=\big[\! -\textstyle\frac{29}{3} ax^3- \textstyle\frac{8}{3}y^3
-(2a^2- \textstyle\frac{5}{2})x^2-9axy+6y^2
+ \textstyle\frac{13}{6} ax+a^2 \big] dx+xy^2dy,
$$
which can be shown to satisfy $I_{\eta_{10}}(h)=0$.
To achieve this, consider the Darboux integral,
$H_{\varepsilon}=\frac{(f_1+\varepsilon g_1)^5}{(f_2+\varepsilon g_2)^4}$,
with
\begin{equation*}
\begin{split}
g_1 =&\ \textstyle\frac{2}{3}x^4+\textstyle\frac{8}{3}ax^3
+\textstyle\frac{4}{3}x^2y+\textstyle\frac{2}{3}x^2
-\textstyle\frac{4}{3}y^2-4ax+4y,\\
g_2 =&\ \textstyle\frac{10}{3}ax^4+\textstyle\frac{5}{3}yx^3
-\textstyle\frac{5}{2}x^3-\textstyle\frac{5}{3}xy^2
-\textstyle\frac{5}{3}ax^2+\textstyle\frac{10}{3}xy
-\textstyle\frac{5}{2}x+\textstyle\frac{10}{3}ay-a,
\end{split}
\end{equation*}
which yields the following system,
\begin{equation*}
\begin{split}
\dot x =&\ M^{-1}H_{0y} - \varepsilon xy^2
+ \varepsilon^2[-\textstyle\frac{2}{9}x^7
+\textstyle\frac{2}{9}ax^6+\textstyle\frac{5}{9}x^5y
-\textstyle\frac{3}{2}x^5-\textstyle\frac{4}{9}ax^4y
-\textstyle\frac{1}{3}x^3y^2\\
        &+\textstyle\frac{17}{9}ax^4+\textstyle\frac{8}{3}x^3y
+\textstyle\frac{2}{9}xy^3-\textstyle(\frac{34}{9}+\frac{16}{9}a^2)x^3
-\textstyle\frac{4}{3}ax^2y-\textstyle\frac{1}{3}xy^2
+\textstyle\frac{2}{9}ax^2\\
        &+\textstyle\frac{7}{3}xy-\textstyle\frac{4}{3}ay^2
-(\textstyle\frac{5}{2}-\textstyle\frac{8}{3}a^2)x
+\textstyle\frac{4}{3}ay-a],\\
\dot y =&-M^{-1}H_{0x} + \varepsilon [ -\textstyle\frac{29}{3} ax^3
- \textstyle\frac{8}{3}y^3 -(2a^2- \textstyle\frac{5}{2})x^2-9axy+6y^2
+ \textstyle\frac{13}{6} ax+a^2 ]\\
&+\varepsilon^2[-\textstyle\frac{4}{9}ax^7-\textstyle\frac{4}{9}x^6y
+(\textstyle\frac{4}{9}a^2+\textstyle\frac{2}{3})x^6
+\textstyle\frac{2}{3}ax^5y+\textstyle\frac{10}{9}x^4y^2
+\textstyle\frac{7}{3}ax^5-2x^4y\\
&-\textstyle\frac{10}{9}ax^3y^2-\textstyle\frac{2}{3}x^2y^3
+(\textstyle\frac{7}{6}-\textstyle\frac{52}{9}a^2)x^4
+\textstyle\frac{13}{9}ax^3y+5x^2y^2+\textstyle\frac{4}{9}y^4
+\textstyle\frac{143}{18}ax^3\\
&-\textstyle\frac{17}{3}x^2y-\textstyle\frac{20}{3}a^2x^2y
-\textstyle\frac{5}{3}axy^2-\textstyle\frac{20}{9}y^3
+(\textstyle\frac{1}{2}+3a^2)x^2-\textstyle\frac{22}{9}axy
+\textstyle\frac{10}{3}y^2\\
&-\textstyle\frac{1}{6}ax-(2-\textstyle\frac{10}{3}a^2)y-a^2].
\end{split}
\end{equation*}
The above system has a center near $C_0$ when $a<-2^{5/4}$
and $|\varepsilon|\ll1$.
Thus, all the Melnikov functions of the above system vanish,
and $M_1(h)=hI_{\eta_{10}}(h)\equiv0$,
implying that $I_{\eta_{10}}(h)=0$.

Next, a direct calculation using $\eta_j$, $1\le j\le 9$, yields
\begin{equation*}
\begin{split}
&\tss\frac{1}{2} (a\eta_1-\eta_4)=(2ax^3-6x^2y+3x^2+2ax)dx+x^3dy\triangleq\bar\eta_4,\\
&\tss\frac{1}{2} (5\eta_1+2\eta_2+6\eta_8)
=(-\textstyle\frac{25}{2}x^3-3ax^2-9xy+\textstyle\frac{1}{2}x+a)dx-ydy\triangleq\bar\eta_2,\\
&\tss\frac{1}{4} (5\eta_1+2\eta_2+2a\eta_3+9\eta_8+2\eta_9)\\
&\quad=(-10x^3-\textstyle\frac{7}{2}xy^2+\textstyle\frac{3}{2}ax^2
-\textstyle\frac{9}{2}xy +\textstyle\frac{5}{8}x+\textstyle\frac{9}{2}ay
+\textstyle\frac{9}{4}a)dx+x^2ydy\triangleq\bar\eta_9,\\
&\tss\frac{1}{12} [2(a^2-10)\eta_1-8\eta_2-14a\eta_3-2a\eta_4-4\eta_7
-21\eta_8-6\eta_9]=[(\textstyle\frac{85}{12}+\textstyle\frac{2}{3}a^2)x^3\\
&-3ax^2y +\textstyle\frac{15}{2}xy^2-4ax^2 -3xy
+(\textstyle\frac{2}{3}a^2-\textstyle\frac{5}{24})x-\textstyle\frac{21}{2}ay
-\textstyle\frac{11}{4}a]dx+y^2dy\triangleq\bar\eta_7,\\
&\tss\frac{1}{2} (3a\eta_1-5\eta_3-\eta_4)
-\eta_5+\eta_6=(ax^3-18x^2-18y-5)dx=a\omega_4 -\omega_5\triangleq\bar\eta_6.
\end{split}
\end{equation*}
Now, suppose we have
$\frac{\omega}{\partial x}=f$ and $\frac{\omega}{\partial y}=g$
for any one-form $\omega=fdx+gdy$. By noticing that
\begin{equation*}
\begin{split}
&\frac{\eta_1}{\partial y}= 1, \ \ \frac{\bar\eta_2}{\partial y}= -y, \ \
\frac{\eta_3}{\partial y}= -x, \ \ \frac{\bar\eta_4}{\partial y}= x^3,  \ \
\frac{\eta_5}{\partial y}= -xy, \ \ \frac{\bar\eta_7}{\partial y}= y^2, \\
&\frac{\eta_8}{\partial y}= x^2, \ \ \frac{\bar\eta_9}{\partial y}= x^2y, \ \
\frac{\eta_{10}}{\partial y}= xy^2, \ \ \frac{\bar\eta_6}{\partial y}= 0, \ \
\frac{\bar\eta_6}{\partial x}\neq 0,
\end{split}
\end{equation*}
we can see that $\eta_1$, $\bar\eta_2$, $\eta_3$, $\bar\eta_4$, $\eta_5$,
$\bar\eta_6$. $\bar\eta_7$, $\eta_8$, $\bar\eta_9$ and $\eta_{10}$
are linearly independent.
Thus, $\eta_{10}$ does not lie in the span of $\eta_j$, $1\le j\le 9$,
and so it follows from $I_{\eta_j}(h)=0$, $j=1,\ldots,10$, that
the dimension of the linear space of integrals
$I_{\omega}(h)$, $\deg(\omega)\le3$ is at most $10$,
Therefore, the independence of integrals $I_{\omega_i}(h)$, $1\le i\le 11$,
proved in \cite{Zoladek1995} does not hold, and there are at most $10$
independent integrals $I_{\omega}(h)$ with $\deg(\omega)\le3$.
The basis can be chosen as $I_{\omega_j}(h)$, $1\le j\le 11$, $j\neq 5$,
since $I_{\bar\eta_6}(h)=\alpha I_{\omega_4}(h)-I_{\omega_5}(h)=0$,
and thus we can remove $I_{\omega_5}(h)$ from the basis
given in \cite{Zoladek1995}.

Regarding $\omega_{12}$, the authors have found a one-form $\bar\omega$,
$\deg(\bar\omega)=3$, by focus value computation such that
the corresponding focus values for
$\widetilde\omega_{12}=\omega_{12}+\bar\omega$ vanish up to a
sufficiently high order. It is reasonable to believe
$I_{\widetilde\omega_{12}}(h)=0$, which implies that using
$\omega_j$, $1\le j\le 12$, can only yield $9$ limit cycles.
However, it is really hard to theoretically prove
$I_{\widetilde\omega_{12}}(h)=0$, yet even this can be done,
we still can not disprove the existence of $11$ small-amplitude limit cycles
for system (\ref{Eq10}), since $\omega_{12}$ is only a special case for
the second-order Melnikov function $M_2(h)$,
and $M_2(h)$ is too complicated to be studied completely.

\section{Cubic Hamiltonian with cubic perturbations}
\label{Hamiltonian}
In order to prove Theorem \ref{TheoB1},
we need some preliminary results for cubic
Hamiltonian given in $(\ref{Eq7})$ with cubic perturbations.
Using the idea and methodology
of \.Zo{\l}\c{a}dek~\cite{Zoladek1995} and~\cite{Iliev2000},
we have the following results summarized in Lemmas 2-5.

Let $\omega_{ij}=x^iy^j dx$ and $\sigma_{ij}=x^iy^j dy$.

\begin{lemma}\label{LemB1}
 For the cubic Hamiltonian given in $(\ref{Eq7})$,
the following identities hold.
\begin{equation*}
   \begin{split}
   (a)\ \sigma_{ij}= \
& \textstyle\frac{1}{j+1}\, d(x^iy^{j+1})
-\textstyle\frac{i}{j+1}\omega_{i-1,j+1};\\
   (b)\  \omega_{ij}= \ &\omega_{i-1,j}
+ \textstyle\frac{j-2i+4}{2j+4}a_1\omega_{i-2,j+2}
-\textstyle\frac{i-2}{j+2}\omega_{i-3,j+2}
-\textstyle\frac{i-2}{j+3}a_2\omega_{i-3,j+3} \\
            & \!\! -x^{i-2}y^{j}dH
+ d \big( \textstyle\frac{1}{j+2}x^{i-2}y^{j+2}
+ \textstyle\frac{a_1}{j+2}x^{i-1}y^{j+2}
+ \textstyle\frac{a_2}{j+3}x^{i-2}y^{j+3} \big),\ i\ge 2;\\
   (c)\ \omega_{0,j}= \
& \textstyle\frac{3j}{a_2(j+1)} \big[
H\omega_{0,j-3}-\textstyle\frac{1}{6}\omega_{1,j-3}
-\textstyle\frac{a_1(j-3)+6j-2}{12(j-1)}
\omega_{0,j-1}-\textstyle\frac{a_1(j+1)}{3(j-1)}\omega_{1,j-1}\\
             &+r_{0,j}(x,y)dH+dR_{0,j}(x,y) \big],\ j\ge 3;\\
   (d)\ \omega_{1,j}= \
& \textstyle\frac{3j}{a_2(j+2)} \big[ H\omega_{1,j-3}
-\textstyle\frac{(j+2)a_1^2}{6(j+1)}\omega_{0,j+1}
+\textstyle\frac{a_2}{6j}\omega_{0,j}
-\textstyle\frac{a_1(5j+3)+6j+2}{12(j-1)}\omega_{1,j-1}\\
   & \!\! - \textstyle\frac{a_1j-3a_1-2}{12(j-1)}\omega_{0,j-1}
-\textstyle\frac{1}{6}
\omega_{1,j-3}
   +r_{1,j}(x,y)dH+dR_{1,j}(x,y) \big],\ j\ge 3;
\end{split}
\end{equation*}
where $r_{i,j}(x,y)$ and $R_{i,j}(x,y)$ are polynomials in $x$ and $y$
with degrees $i+j-2$ and $i+j+1$, respectively.
\end{lemma}

\begin{proof}\
A direct calculation using integration by parts results in the formula (a).
From the Hamiltonian, we have the equation
$\frac{1}{3}x^3=\frac{1}{2}(x^2+y^2)+\frac{1}{2}a_1xy^2
+\frac{1}{3}a_2y^3-H$, giving the relation,
\begin{equation*}
x^2dx=xdx+ydy+\frac{a_1}{2}y^2dx+a_1xydy+a_2y^2dy-dH,
\end{equation*}
which in turn yields
\begin{equation}\label{Eq16}
\omega_{i,j} = \omega_{i-1,j}+\sigma_{i-2,j+1}
+ \frac{a_1}{2}\omega_{i-2,j+2} + a_1\sigma_{i-1,j+1}
+ a_2\sigma_{i-2,j+2} - x^{i-2}y^jdH,\ i\ge 2.
\end{equation}
Then, combining (\ref{Eq16}) with the formula (a) we obtain the formula (b).

Similarly, the equation,
$\frac{1}{3}a_2y^3=H-\frac{1}{2}(x^2+y^2)+\frac{1}{3}x^3
-\frac{1}{2}a_1xy^2$, generates
\begin{equation}\label{Eq17}
\frac{1}{3}a_2\omega_{i,j}=H\omega_{i,j-3}-\frac{1}{2}\omega_{i+2,j-3}
-\frac{1}{2}\omega_{i,j-1}+\frac{1}{3}\omega_{i+3,j-3}
-\frac{1}{2}a_1\omega_{i+1,j-1},\ j\ge 3.
\end{equation}
Finally, the formulas (c) and (d) follow the formula (b) and (\ref{Eq17}).
\end{proof}

From Lemma \ref{LemB1}, we know that any polynomial
one-form $\omega$, $\deg(\omega)=m$, can be expressed in the form of
$$
\omega=r(x,y)dH+dR(x,y)
+\sum_{i=0,1}\sum_{j=0}^{m-i}\alpha_{i,j}\omega_{i,j}.
$$
The next lemma shows that there also exist some relationships among
the one-forms $\omega_{i,j}$,
$i=0,1$.

\begin{lemma}\label{LemB2}
 For any non-negative integer, $m\bmod 3\neq 2$,
there exist $\beta_{i,j,m}$, $\widetilde r_m(x,y) $ and
$\widetilde R_m(x,y)$ satisfying the following identity,
\begin{equation}\label{Eq18}
\sum_{i=0,1}\sum_{j=0}^{m-i}\beta_{i,j,m}\omega_{i,j}
=\widetilde r_m(x,y)dH+d\widetilde R_m(x,y),
\end{equation}
where $\widetilde R_m(x,y)$ and $\widetilde r_m(x,y)$ are polynomials
in $x$ and $y$ of degrees $m+1$ and $m-1$, respectively; and $\beta_{i,j,m}$
are polynomials in $a_{1}$ and $a_{2}$, with $\beta_{0,0,0}=\beta_{1,0,1}=1$,
$\beta_{0,1,1}=0$, and
\begin{equation}\label{Eq19}
\begin{split}
\beta_{0,m+3,m+3}=&\frac{m+4}{3(m+3)}(a_2\beta_{0,m,m}
+\frac{a_1^2}{2}\beta_{1,m-1,m}),\\
\beta_{1,m+2,m+3}=&\frac{m+4}{3(m+2)}(a_1\beta_{0,m,m}+a_2\beta_{1,m-1,m}),
\end{split}
\end{equation}
if $\beta_{1,-1,0}$ is defined as $\beta_{1,-1,0}=0$.
\end{lemma}

\begin{proof}
We use the method of mathematical induction to prove this lemma. It is easy
to see that the conclusion is true for $m=0,1$. Now, suppose (\ref{Eq18})
holds for $m\bmod 3\neq 2$. Then, we prove that (\ref{Eq18}) also holds
for $m+3$. Multiplying (\ref{Eq18}) by $H$ on both sides yields
\begin{equation}\label{Eq20}
\sum_{i=0,1}\sum_{j=0}^{m-i}\beta_{i,j,m}H\omega_{i,j}
=H\widetilde r_mdH+Hd\widetilde R_m.
\end{equation}
The right-hand side of (\ref{Eq20}) can be rewritten as
\begin{equation}\label{Eq21}
H\widetilde r_mdH+Hd\widetilde R_m
=(H\widetilde r_m-\widetilde R_m)dH+d(H\widetilde R_m).
\end{equation}
For the left-hand side of (\ref{Eq20}), it follows from the formulas (c)
and (d) in Lemma \ref{LemB1} that
\begin{equation}\label{Eq22}
\begin{split}
&H\omega_{i,j}=\xi_{i,j+3}+\eta_{i,j+3},\ i+j<m,\\
&H\omega_{0,m}=\frac{a_2(m+4)}{3(m+3)}\omega_{0,m+3}
+\frac{a_1(m+4)}{3(m+2)}\omega_{1,m+2}+\eta_{0,m+3},\\
&H\omega_{1,m-1}=\frac{a_1^2(m+4)}{6(m+3)}\omega_{0,m+3}
+\frac{a_2(m+4)}{3(m+2)}\omega_{1,m+2}+\eta_{1,m+2},\ m>0,
\end{split}
\end{equation}
where $\eta_{i,j}=r_{i,j}dH+dR_{i,j}$, and $\xi_{i,j}$ is a one-form with
$\deg(\xi_{i,j})\le i+j$. Then, substituting (\ref{Eq22}) into the
left-hand side of (\ref{Eq20}) yields
\begin{equation}\label{Eq23}
\begin{split}
\sum_{i=0,1}\sum_{j=0}^{m-i}\beta_{i,j}H\omega_{i,j}
=&\frac{m+4}{3(m+3)}(a_2\beta_{0,m,m}
+\frac{a_1^2}{2}\beta_{1,m-1,m})\omega_{0,m+3}\\
&+\frac{m+4}{3(m+2)}(a_1\beta_{0,m,m}+a_2\beta_{1,m-1,m})\omega_{1,m+2}\\
&+\sum_{i=0,1}\sum_{j=0}^{m-i}\beta_{i,j}(\xi_{i,j+3}+\eta_{i,j+3}).
\end{split}
\end{equation}
Finally, combining (\ref{Eq23}) with (\ref{Eq20}) and (\ref{Eq21}) shows
that the conclusion is also true for $m+3$.

The proof of the lemma is complete.
\end{proof}

Noting that $\beta_{0,0,0}=\beta_{1,0,1}=1$,
$\beta_{1,-1,0}=\beta_{0,1,1}=0$, we know from (\ref{Eq19}) that
$\beta_{k,m-k,m}$ in Lemma \ref{LemB2} are polynomials in $a_1$ and $a_2$
with positive coefficients for $m\bmod 3=k$, $k<2$. Thus,
it follows from (\ref{Eq18}) that $\omega_{k,m-k}$, $m\bmod 3=k<2$,
can be expressed in terms of other one-forms $\omega_{i,j}$, $i+j\le m$
and $r_mdH+dR_m$. This gives the following lemma.

%\vspace{0.10in}
%\noindent
\begin{lemma}\label{LemB3}
Any polynomial one-form $\omega$ of degree $m$ can
be expressed as
\begin{equation}\label{Eq24}
\omega=r(x,y)dH+dR(x,y)+\sum_{i=0,1}\sum_{j\bmod 3\neq 0}^{1\le j\le m-i}
\alpha_{ij}\omega_{ij},
\end{equation}
where $R(x,y)$ and $r(x,y)$ are polynomials in $x$ and $y$
of degrees $m+1$ and $m-1$, respectively.
\end{lemma}

\vspace{0.10in}
Now, we use (\ref{Eq24}) to obtain
\begin{equation}\label{Eq25}
M(h)=\oint_{\Gamma_h}\omega=\sum_{i=0,1}\sum_{j\bmod 3\neq 0}^{1\le j\le m-i}
\alpha_{ij}\oint_{\Gamma_h}\omega_{ij},
\end{equation}
which implies that any Melnikov function $M(h)=\oint_{\Gamma_h}\omega$,
$\deg(\omega)=m$, can be expressed as a linear combination of integrals
$I_{ij}(h)=\oint_{\Gamma_h}\omega_{ij}$, $i=0,1$, $j \bmod 3\neq0$.
A reasonable expectation is that the integrals $I_{i,j}(h)$ form a basis
for the linear space of Melnikov functions $M(h)=\oint_{\Gamma_h}\omega$.
Actually, it will be seen in the next section that the space of Melnikov
functions $M(h)$ could be Chebyshev with accuracy at least 2.
So the number of limit cycles in system (\ref{Eq6}) is not determined
by the number of elements in the basis.
Further, the coefficients $\alpha_{i,j}$ in (\ref{Eq25}) could become
very complicated when $M(h)$ is a higher-order Melnikov function of
system (\ref{Eq6}). In this case, it is really not easy to prove the
independence of $\alpha_{i,j}$s, which is the second big obstacle in the
use of the independence of the integrals $I_{i,j}(h)$ to determine the
number of limit cycles.

To overcome the above mentioned difficulty, we turn to an alternative,
which decreases the complexity in computing $M(h)$ by (\ref{Eq24}),
but it still does not solve the problem of independence of basis.
Let $\omega_j=Q_j(x,y)dx-P_j(x,y)dy$. Then, for higher-order Melnikov
functions of system (\ref{Eq6}), we have the following result.

%\vspace{0.10in}
%\noindent
\begin{lemma}\label{LemB4} (cf. \cite{Iliev2000,Francoise1996})
Let $(\ref{Eq9})$ hold. Assume that in system (\ref{Eq6})
for some $k\ge2$, Melnikov function
$M_m(h)=\oint_{\Gamma_h}\Phi_m\equiv0$, $1\le m\le k-1$, and $\Phi_m$
can be expressed as
\begin{equation}\label{Eq26}
\Phi_m=r_mdH+dR_m.
\end{equation}
Then,
\begin{equation}\label{Eq27}
\begin{split}
M_k(h) &= \oint_{\Gamma_h} \Big(\omega_k+\sum_{i+j=k}r_i\omega_j \Big),\\
r_mdH+dR_m &= \omega_m +\sum_{i+j=m}r_i\omega_j,\ 1\le m\le k-1.
\end{split}
\end{equation}
\end{lemma}

\begin{proof}
We prove this lemma by using the method of mathematical induction.
First, write system (\ref{Eq6}) in the Pfaffian form,
\begin{equation}\label{Eq28}
dH-\varepsilon \omega_1 -\varepsilon^2\omega_2 - \cdots=0.
\end{equation}
Multiplying (\ref{Eq28})
by $1+\varepsilon r_1+\ldots+ \varepsilon^{k-1}r_{k-1}$ and combing
the like terms yield
\begin{equation*}
\begin{split}
dH&+\varepsilon(r_1dH-\omega_1)+\varepsilon^2(r_2dH-r_1\omega_1-\omega_2)
+\cdots\\
&+\varepsilon^k(-r_{k-1}\omega_1-\cdots-r_1\omega_{k-1}
-\omega_k)+O(\varepsilon^{k+1})=0,
\end{split}
\end{equation*}
which, by using (\ref{Eq27}), can be written as
\begin{equation*}
dH-\varepsilon dR_1-\cdots-\varepsilon^{k-1}dR_{k-1}-
\varepsilon^k(r_{k-1}\omega_1+\cdots+r_1\omega_{k-1}+\omega_k)
+O(\varepsilon^{k+1})=0.
\end{equation*}
Then, we integrate the above equation along the phase curve $\gamma$
from point $A$ to point $B$, which are used to define the first return map.
Note that
$$
d(h,\varepsilon)=\int_{\gamma}dH=H(B)-H(A)=O(|A-B|)
$$
and
$$
\Big|\int_{\gamma}(\varepsilon dR_1+\varepsilon^2dR_2+\cdots
+\varepsilon^{k-1}dR_{k-1})\Big|=\varepsilon \, O(|A-B|).
$$
In addition, it follows from (\ref{Eq8}) that
$d(h,\varepsilon)=O(\varepsilon^k)$. Therefore,
$|A-B|=O(\varepsilon^k)$ and we finally obtain
\begin{equation*}
d(h,\varepsilon)=\varepsilon^k\int_\gamma(r_{k-1}\omega_1
+\cdots+r_1\omega_{k-1}+\omega_k)+O(\varepsilon^{k+1}),
\end{equation*}
which yields
$$
M_k(h)=\oint_{\Gamma_h} \Big(\omega_k+\sum_{i+j=k}r_i\omega_j\Big).
$$
The proof is finished.
\end{proof}

\begin{rmk}
For the generic (system) parameters $(a_1,a_2)\in\Omega$,
system (\ref{Eq6}) satisfies Fran\c{c}oise's $\ast$-property
\cite{Jebrane2007}: for any polynomial one-form $\omega$,
if $\oint_{\Gamma_h}\omega\equiv0$, then $\omega=rdH+dR$
for some polynomials $r$ and $R$.
So the only condition which is needed in Lemma \ref{LemB4} is $M_m(h)\equiv0$
when generic Hamiltonians are considered.
\end{rmk}

\begin{rmk}
For some cubic Hamiltonians, the Fran\c{c}oise's $\ast$-property does
not hold (see \cite{Uribe2006}).
In other words, in such systems we could have polynomial one-forms $\omega$
satisfying $\oint_{\Gamma_h}\omega\equiv0$,
but $\omega$ can not be expressed in the form of $\omega=rdH+dR$,
where $r$ and $R$ are some polynomials.
Therefore, it is required that $\Phi_m$ should not contain
such ``bad'' one-forms for Melnikov function
$M_m(h)=\oint_{\Gamma_h}\Phi_m\equiv0$ in Lemma \ref{LemB4}.
\end{rmk}

\section{Proof of Theorem \ref{TheoB1}}
\label{proof_Thm1}

Now with the results obtained in the previous section,
we are ready to prove Theorem \ref{TheoB1}.

\begin{proof}
We return to system (\ref{Eq6}) with $P(x,y)$ and $Q(x,y)$ defined in
(\ref{Eq9}), and want to use higher-order Melnikov functions to prove
the existence of $10$ small-amplitude limit cycles around the origin.

Due to the difficulty in the proof of independence of basis, we use the
computation of focus values to prove the theorem. However, the computation
becomes very demanding or almost impossible for computing higher-order
focus values if all the coefficients are retained in the computation,
and in fact many terms are not necessarily needed. Thus, before computing
the focus values of system (\ref{Eq6}), without loss of generality, we want
to simplify this system by choosing a group of coefficients $a_{ijm}$,
$b_{ijm}$ in the polynomials $P(x,y)$ and $Q(x,y)$, which does not reduce
the number of limit cycles bifurcating from the origin.

In the following, we shall show how to choose a group of coefficients
which are necessary for the first non-vanishing Melnikov function $M_k(h)$
in (\ref{Eq8}). Based on the results presented in the previous section
(in particular, Lemmas \ref{LemB1}, \ref{LemB3} and \ref{LemB4}), we provide an algorithm as follows.

Consider $M_1(h)$ in system (\ref{Eq6}), we know
$M_1(h)=\oint_{\Gamma_h}\omega_1$. Using Lemma \ref{LemB3}, we have
\begin{equation}\label{Eq29}
\omega_1=Q_1dx-P_1dy=\sum_{i=0}^1\sum_{j=1}^2\alpha_{ij1}x^iy^jdx+r_1dH+dR_1,
\end{equation}
with $r_1=-(b_{211}+3a_{301})y$. Then,
\begin{equation*}
M_1(h)=\oint_{\Gamma_h} \big( \alpha_{011}ydx+\alpha_{111}xydx
+\alpha_{021}y^2dx+\alpha_{121}xy^2dx \big).
\end{equation*}
It is seen that $M_1(h)$ depends on $\alpha_{ij1}$, $i=0,1$, $j=1,2$.
So only four coefficients in the polynomials $P_1(x,y)$ and $Q_1(x,y)$
are needed in order to keep $\alpha_{ij1}$, $i=0,1$, $j=1,2$
being independent without decreasing the number of zeros of $M_1(h)$.
We choose these four coefficients as $b_{ij1}$, $i=0,1$, $j=1,2$.
(Certainly, the choice is not unique.) Then, we have polynomials
\begin{equation}\label{Eq30}
P_1(x,y)=0,\quad Q_1(x,y)=b_{011}x+b_{111}xy+b_{021}y^2+b_{121}xy^2.
\end{equation}

Next, let us consider $M_2(h)$ when
$M_1(h)=\oint_{\Gamma_h}r_1dH+dR_1\equiv0$, i.e., all $\alpha_{ij1}=0$
in (\ref{Eq29}). Lemma \ref{LemB4} gives $M_2(h)=\oint_{\Gamma_h}
\widetilde\omega_2$,
where $\widetilde\omega_2=\omega_2+r_1\omega_1$.
Thus, by using Lemma \ref{LemB3}, we obtain
\begin{equation*}
\widetilde\omega_2=\sum_{i=0}^1\sum_{j=1}^2\alpha_{ij2}x^iy^jdx
+\alpha_{042}y^4dx+r_2dH+dR_2,
\end{equation*}
which shows that $M_2(h)$ depends on $\alpha_{ij2}$, $i=0,1$, $j=1,2$
and $\alpha_{042}$. Obviously, the coefficient $\alpha_{042}$ is derived
from $r_1 \omega_1$ by Lemma \ref{LemB3} because the one-form $y^4dx$
of degree 4 comes from $r_1 \omega_1$. For $\varepsilon$-order perturbations,
$b_{ij1}$, $i=0,1$ $j=1,2$ are needed to get all $\alpha_{ij1}=0$
in (\ref{Eq29}). For $r_1$ we may simply take $b_{211}=1$ and $a_{301}=0$,
yielding $r_1=-y$. We also see that the one-form $y^4dx$ can be derived
from $x^3ydx$ by using the formula (b) in Lemma \ref{LemB1}.
Hence, we may choose $b_{301}$ for $\alpha_{042}$
so that $b_{301}x^3ydx$ could appear in $r_1\omega_1$. For $\alpha_{ij2}$,
$i=0,1$, $j=1,2$, by an argument similar to that for $M_1(h)$, we choose
$b_{012}$, $b_{112}$, $b_{022}$ and $b_{122}$. Hence, we obtain the
following polynomials,
\begin{equation}\label{Eq31}
\begin{split}
\!\! &P_1(x,y)=0,\quad Q_1(x,y)=b_{011}x \!+\! b_{111}xy \!+\! b_{021}y^2
\!+\! b_{121}xy^2 \!+\! b_{301}x^3 \!+\! x^2y,\\
\!\! &P_2(x,y)=0,\quad Q_2(x,y)=b_{012}x+b_{112}xy+b_{022}y^2+b_{122}xy^2.
\end{split}
\end{equation}

Following the above procedure, we can choose the coefficients for $M_3(h)$,
and so on. In the following, we list the polynomials for $M_k(h)$ up to $k=5$
(the detailed arguments are omitted here for brevity):
\begin{equation}
\begin{split}\label{Eq32}
P_j(x,y) =\ & a_{21j}x^2y+a_{12j}xy^2,\ j=1,2,3,\quad P_4(x,y) =P_5(x,y) = 0,\\
Q_1(x,y) =\ & b_{011}y \!+\! b_{111}xy \!+\! b_{021}y^2 \!+\! b_{121}xy^2
\!+\! b_{301}x^3 \!+\! b_{031}y^3 \!+\! b_{211}x^2y,\\
Q_2(x,y) =\ & b_{012}y+b_{112}xy+b_{022}y^2+b_{122}xy^2+b_{302}x^3+b_{032}y^3,\\
Q_3(x,y) =\ & b_{013}y+b_{113}xy+b_{023}y^2+b_{123}xy^2+b_{303}x^3,\\
Q_4(x,y) =\ & b_{014}y+b_{114}xy+b_{024}y^2+b_{124}xy^2+b_{304}x^3,\\
Q_5(x,y) =\ & b_{015}y+b_{115}xy+b_{025}y^2+b_{125}xy^2.
\end{split}
\end{equation}
Here, the difficult part is to compute the functions
$r_i, \, i=1,2,3,4$ in $\widetilde\omega_i$. Proving the independence of
the basis for each $k$ is even more difficulty.
Thus, we turn to focus value computation which can be easily used to
show the independence of the basis (i.e. the focus values).

Having determined the coefficients we need in $P$ and $Q$ of
system (\ref{Eq6}), we now use the computation of focus values to prove
the existence of $10$ small-amplitude limit cycles. We compute the focus
values up to $\varepsilon^5$ order as follows:
\begin{equation}\label{Eq33}
V = \sum_{i=0}^5 \varepsilon^i V_i, \quad
{\rm where} \quad
V_i = \{v_{i0}, v_{i1}, v_{i2}, \cdots \}.
\end{equation}
We call $v_{ij}$ the $j$th $\varepsilon^i$-order focus value of system
(\ref{Eq6}), and note that $v_{0j}=0, \, j=0,1,2, \dots$ since at
$\varepsilon=0$ system (\ref{Eq6}) is a Hamiltonian system. The computation
of $V_i$ is equivalent to the computation of $i$th-order Melnikov function
$M_i(h)$. But the computation of focus values is much easier than that
of the higher-order Melnikov functions. The disadvantage of the focus value
computation is that conditions obtained from the first few focus values
are hard to be used to prove vanishing of an infinite number of focus values.
But this can be easily verified by the above formulas $\widetilde\omega_i$.

The focus values $v_{ij}$ can be obtained by using many different
symbolic programs (e.g., the Maple program developed in~\cite{Yu1998}).
Firstly, note that $v_{i0} = \frac{1}{2}\, b_{01i}, \, i=1,2, \dots$.
In order to execute the Maple program, set $b_{01i} = 0, \, i=1,2, \dots$.
In addition, set $b_{211} = 1$. Now, we start from $V_1$ and obtain
$$
v_{11} = \textstyle\frac{1}{8} ( a_{121}+3 b_{031}+ b_{111}-
\textstyle\frac{1}{2} a_1 b_{111}-2 a_2 b_{021} +1 ) .
$$
Setting $v_{11}=0$ yields
$  b_{031} = ( \frac{1}{2} a_1 b_{111}+2 a_2 b_{021} -a_{121}- b_{111} -1)/3.
$
Further, setting $v_{12} = 0$ results in
$$
b_{121} = a_1  b_{021}  - a_{211}
+ \frac{1}{4 a_2 (5 a_1 -2)}
(3 a_1^2 +20 a_2^2 +4 a_1 -20) (b_{111} +1).
$$
Then, we have
\begin{equation*}
\begin{split}
v_{13} =&\, \frac{35}{3072 (5 a_1-2)} (b_{111}+1)
(a_1^3-3 a_1^2+4-4 a_2^2) F_{11}, \\
v_{14} =&\, \frac{-7}{73728 (+5 a_1-2)}
(b_{111}+1) (a_1^3-3 a_1^2+4-4 a_2^2) F_{12}, \\
v_{15} =&\, \frac{-7}{84934656 (+5 a_1-2)}
(b_{111}+1) (a_1^3-3 a_1^2+4-4 a_2^2) F_{13},
\end{split}
\end{equation*}
where
$$
\begin{array}{ll}
F_{11} = 3 a_1^2 +12 a_1 - 4 - 4 a_2^2, \\ [0.5ex]
F_{12} = 27 a_1^4 \!-\! 90 a_1^3 \!-\! 1308 a_1^2 \!+\! 1608 a_1 \!-\!256
\!+\! ( 420 a_1^2 \!+\! 1608 a_1 \!-\! 1376 \!-\! 256 a_2^2 ) a_2^2 , \\ [0.5ex]
F_{13} = 19683 a_1^6 \!+\! 343116 a_1^5 \!-\! 124524 a_1^4 \!-\!
6168672 a_1^3 \!+\! 7612368 a_1^2 \!+\! 1585344 a_1 \\ [0.5ex]
\qquad \quad
-1071424 \!+\! 4 [ 3(140715 a_1^4 \!+\! 622536 a_1^3 \!+\!  39880 a_1^2
\!-\! 1689568 a_1 \!+\! 421808)\\[0.5ex]
\qquad \quad - (404508 a_1^2-396336 a_1+267856 a_2^2-1265424) a_2^2 ] a_2^2 .
\end{array}
$$
It is easy to see that setting $b_{111}=-1$ results in
$v_{13} = v_{14} = v_{15} = \cdots = 0$, as discussed above.
In order to obtain maximal number of small-amplitude limit cycles
bifurcating from the origin, we have to use the coefficients
$a_1$ and $a_2$ to solve $F_{11}=F_{12}= 0$ (i.e., $v_{13} = v_{14} = 0$).
If the solution of $F_{11}=F_{12} = 0$ yields $F_{13} \ne 0$, i.e.,
we have parameter values such that
$v_{10} = v_{11} = \cdots = v_{14} = 0$, but $v_{15} \ne 0$,
then we obtain $5$ small-amplitude limit cycles by properly perturbing
$b_{011}$, $ b_{031}, \, b_{021}, \, a_1 $ and $a_2$, respectively. In fact,
by using the Groebner basis reduction procedure, we can reduce $F_{12}$
and $F_{13}$ to
\begin{equation*}
\begin{split}
\tilde{F}_{12} =&\, \left. F_{12} \right|_{F_{11}=0}
= 18(a_1+2) ( 11 a_1^3 + 46 a_1^2 - 84 a_1 + 24), \\
\tilde{F}_{13} =&\, \left. F_{13} \right|_{F_{11}= \tilde{F}_{12} = 0}
= - \tss\frac{179712}{121} (a_1 + 2) (3073 a_1^2 - 5272 a_1 + 1500) \ne 0.
\end{split}
\end{equation*}
In fact, solving the system of two equations,
$F_{11} = F_{12} = 0 $ (or $F_{11}= \tilde{F}_{12} = 0$) we obtain the
solutions for $a_1$ as follows:
\begin{equation}
\begin{array}{ll}
a_1 = a_{11}^i, \quad a_2 = a_{21}^i = \pm \frac{1}{2} \sqrt{ 3 (a_{11}^i)^2
+ 12 a_{11}^i - 4}, \ \ i = 1,2,3, \ \ {\rm for \ which} \\ [1.0ex]
a_{11}^1=-5.61185383 \cdots, \quad
a_{11}^2 = 0.36507058 \cdots, \quad a_{11}^3 = 1.06496506 \cdots ,
\label{Eq_soln1}
\end{array}
\end{equation}
where the second number `1' in the subscripts of $a_{11}^i$ and $a_{21}^i$
denotes the solutions corresponding to the first-order Melnikov function,
i.e, $k=1$.
Note that $a_1=-2$ is not a solution of $F_{11}=0$. Further, we obtain
%$$
%\det \! \left[ \begin{array}{cc}
%\displaystyle\frac{\partial F_{11}}{\partial a_1} &
%\displaystyle\frac{\partial F_{11}}{\partial a_2} \\
%\displaystyle\frac{\partial F_{12}}{\partial a_1} &
%\displaystyle\frac{\partial F_{12}}{\partial a_2} \end{array}
%\right]_{F_{11}= \tilde{F}_{12} = 0}
%=  576 a_2 (a_1+1) (11 a_1^2+40 a_1-36) \ne 0,
%$$
$$
\det \left[\frac{\partial (F_{11},F_{12})}{\partial (a_1,a_2)}\right]_{F_{11}
= \tilde{F}_{12} = 0}
=  576 a_2 (a_1+1) (11 a_1^2+40 a_1-36) \ne 0,
$$
since none of the factors in the above equations
are included in $F_{11}$ and $\tilde{F}_{12}$.
%verified by directly substituting the solutions given in
%(\ref{Eq_soln1}) into the above determinant.

Summarizing the above results we can conclude that based on the
$\varepsilon^1$-order focus values (equivalently based on the first-order
Melnikov function $M_1(h)$) we obtain $5$ small-amplitude limit
cycles around the origin.

Now let $ b_{111} = -1$, then $b_{121} = a_1 b_{021} - 1$ and
$ b_{031}=-\frac{1}{3} ( a_{121} + \frac{1}{2} a_1 -2 a_2 b_{021} )$,
under which all $\varepsilon^1$-order focus values vanish, or
equivalently, the first-order Melnikov function $M_1(h) \equiv 0$.
Note here that $a_1$ and $a_2$ are not used in
making $M_1(h) \equiv 0$.
Then, one uses the $\varepsilon^2$-order focus values to solve the
polynomial equations $v_{21}=v_{22} = v_{23} = 0 $, yielding the solutions
for $ b_{032}, \, b_{122} $ and $b_{112}$. Under these solutions,
we further obtain
$$
\begin{array}{ll}
v_{24} = - &\!\!\!\!\! \tss\frac{1}{36864(3 a_1^2+12a_1-4-4a_2^2)}
\ F_{20} F_{21}, \\ [1.0ex]
v_{25} =   &\!\!\!\!\! \tss\frac{1}{31850496 (3 a_1^2+12a_1-4-4a_2^2)}
\ F_{20} F_{22}, \\ [1.0ex]
v_{26} = &\!\!\!\!\! \tss\frac{11}{107297229312 (3 a_1^2+12a_1-4-4a_2^2)}
\ F_{20} F_{23},
\end{array}
$$
for which we have applied the Groebner basis reduction procedure to obtain
$$
\begin{array}{ll}
F_{20} = [ 2 (3 a_1^3 \!-\! 4 a_2^2) b_{021}
- 3 ( a_1^3 \!-\! 4 a_2^2)  b_{301}
\!-\! 6 a_1^2 a_{211}
\!+\! 4 a_2 a_{121}
\!-\! 4 a_1 a_2 b_{211} ] b_{211}  \\ [0.5ex]
\qquad \ + 12 a_1 ( a_1 a_{121} - a_2 a_{211}) b_{021},  \\[0.5ex]
F_{21} =
81 a_1^4 \!-\! 648 a_1^3 \!-\! 648 a_1^2 \!+\! 1632 a_1 \!-\! 880
\!-\! (504 a_1^2 \!-\! 1632 a_1 \!-\! 1696 \!+\! 880 a_2^2) a_2^2 , \\[0.5ex]
\tilde{F}_{22} = \left. F_{22} \right|_{F_{21}=0} \\[0.5ex]
\hspace{0.285in} = 1408 [ 243 a_1^3 \!-\! 522 a_1^2 \!+\! 5172 a_1
\!+\! 6664 \!+\! (1053 a_1^2 \!-\! 2424 a_1 \!-\! 5572 \!+\! 1300 a_2^2) a_2^2
] a_2^2\\ [0.5ex]
\qquad \ - 50688 ( 63 a_1^3 + 56 a_1^2 - 148 a_1 + 80), \\[0.5ex]
\tilde{F}_{23} = \left. F_{23} \right|_{F_{21}=\tilde{F}_{22}=0} \\[0.5ex]
\hspace{0.285in} =
72 (675121644 a_1^3+475639745 a_1^2-1491227668 a_1+849702020) \\[0.5ex]
\qquad \
+ \{ 3893155245 a_1^3 +22056197796 a_1^2 -131201934348 a_1
        -117343356608 \\[0.5ex]
\qquad \quad
+20 [ 303274623 a_1 +3083354476
          -26 (55458 a_1 -130879) a_2^2 ] a_2^2 \} a_2^2 \ne 0,
\end{array}
$$
Similarly, we obtain the following solutions satisfying $F_{21} =
\tilde{F}_{22} = 0$:
\begin{equation}
\begin{array}{ll}
a_1 = a_{12}^i, \quad i=1,2, \dots, 7, \\[0.5ex]
a_2 = a_{22}^i \!=
\! \sqrt{
\textstyle\frac{10179 a_1^6
-81864 a_1^5
-179172 a_1^4
+204992 a_1^3
-32496 a_1^2
-124032 a_1
+66880}
{4 ( 5109 a_1^4+12076 a_1^3-75936 a_1^2-167664 a_1+48944 )}
}, \ (a_1 = a_{12}^i), \\[0.5ex]
{\rm where} \\[0.5ex]
a_{12}^1 = -2.43192492 \cdots, \ \
a_{12}^2 =  0.12148877 \cdots, \ \
a_{12}^3 =  0.23963547 \cdots, \\[0.5ex]
a_{12}^4 = \hspace{0.135in}  0.89471272 \cdots, \ \
a_{12}^5 =  1.60031174 \cdots, \ \
a_{12}^6 =  7.33752703 \cdots, \\[0.5ex]
a_{12}^7 = \hspace{0.06in} 10.40950390 \cdots.
\end{array}
\label{Eq_soln2}
\end{equation}

In addition, we can show that for the above solutions the following
determinant is non-zero,
\begin{equation*}
\begin{split}
&\det \left[\frac{\partial (F_{21},F_{22})}{\partial (a_1,a_2)}\right]_{F_{21}
= \tilde{F}_{22} = 0}\\
= &\, \frac{360448}{351} a_2 \{
36 (1571445 a_1^3+860083 a_1^2-3207848 a_1+1911580) \\
&\,\qquad \qquad + a_2^2 [ 4977612 a_1^3+24045705 a_1^2
-138196596 a_1-132836684 \\
&\, \qquad \qquad\qquad + 20 a_2^2 (-119877 a_1+2945227
+ 169 a_2^2 (459 a_1+1799) ) ] \} \ne 0.
\end{split}
\end{equation*}
The above results show that we have parameter values such that
$\, v_{20}= v_{21} = \cdots = v_{25} = 0$, but $v_{26} \ne 0$.
Then, taking proper perturbations on the coefficients
$b_{012}$, $ b_{032}$, $ b_{122}$, $ b_{112}$, $ a_1 $ and $a_2$
yields $6$ small-amplitude limit cycles around the origin of system
(\ref{Eq6}) when the $\varepsilon^2$-order focus values
(or the second-order Melnikov function $M_2(h)$) are used.

In order to get more limit cycles, we let $F_{20} = 0$ and solve
this equation for $b_{301}$, yielding all the $\varepsilon^2$-order focus
values $v_{2j}=0$. Under these conditions, we then use the
$\varepsilon^3$-order focus values $v_{3j}$ to determine the number of
small-amplitude limit cycles. Similarly, we may linearly solve the
polynomial equations $v_{31} = v_{32} = v_{33} = v_{34} = 0 $ for the
coefficients $b_{023}$, $ b_{123}$, $ b_{113}$ and $b_{302}$. After this,
no coefficients can be solved linearly. So we solve $a_{211}$ from
the equation, $v_{35}=0$, which is quadratic about $a_{211}$,
to obtain two solutions $a_{211}^\pm$.
We choose $a_{211} = a_{211}^+$ and then $v_{36}$, $v_{37}$ and
$v_{38}$ are simplified to
$$
v_{36} = -624\, F_{30}\, F_{31}, \quad
v_{37} = -1248\, F_{30}\, F_{32}, \quad
v_{38} = -208\, F_{30}\, F_{33},
$$
where $F_{30}$ is a lengthy irrational function, and
we further apply the Groebner reduction procedure to $F_{32}$
and $F_{33}$ to obtain
\begin{equation*}
\begin{split}
F_{31} \!= & \,
405 a_1^4 \!\!+\! 6264 a_1^3 \!\!+\! 6264 a_1^2 \!-\! 5664 a_1 \!+\! 1360
     \!-\! 8 (99 a_1^2 \!+\! 708 a_1 \!+\! 524 \!-\! 170 a_2^2) a_2^2 ,
\\
\tilde{F}_{32} \!= &  \left. F_{32} \right|_{F_{31}=0} \\
\!= &  \,
4 ( 261117 a_1^3 + 307422 a_1^2 - 260532 a_1 + 60680)
-[ 9 (1035 a_1^3 + 13266 a_1^2 \\
& + 111492 a_1 + 84376
+ 5 (513 a_1^2 -4824 a_1 -57156 +2660 a_2^2) a_2^2 ] a_2^2, \\
\tilde{F}_{33} \!= &
\left. F_{32} \right|_{F_{31}=\tilde{F}_{32}=0} \\
\!= &  \, 4 (152348063679 a_1^3+175217936814 a_1^2-151386504684 a_1
\\
& +35757329960)
+ \{ 7428338685 a_1^3 -38896637238 a_1^2 \\
& -568264627476 a_1 -439876872808
- 20 [ 714254595 a_1 -6998804702 \\
& - 380 (11970 a_1 + 132193 ) a_2^2 ] a_2^2 \} a_2^2 \ne 0.
\end{split}
\end{equation*}
Solving $F_{31}=\tilde{F}_{32}=0$ yields
\begin{equation}
\begin{array}{ll}
a_1 = a_{13} = 0.01871627 \cdots, \\ [0.5ex]
a_2 = a_{23} = \pm \sqrt{\frac{99 a_{13}^2 + 708 a_{13} + 524
- 12 \sqrt{1104 + 8496 a_{13} +504 a_{13}^2 - 2724 a_{13}^3
- 171 a_{13}^4 } }{340}
}.
\end{array}
\label{Eq_soln3}
\end{equation}

Further, we have
$
%\begin{array}{rl}
%& \!\!\!
\det \left[\partial (F_{31},F_{32})/\partial (a_1,a_2) \right]
%\det \! \left[ \begin{array}{cc}
%\frac{\partial F_{31}}{\partial a_1} &
%\frac{\partial F_{31}}{\partial a_2} \\
%\frac{\partial F_{32}}{\partial a_1} &
%\frac{\partial F_{32}}{\partial a_2} \end{array}
%\right]_{(a_1,a_2)=(a_{13},a_{23}) }
= -0.1124026367 \cdots \times 10^{10} \ne 0
%  \\[2.5ex]
%= & \!\!\! \frac{2228224}{475}\, a_2 \{
%238627593 a_1^3+254229138 a_1^2-235287828 a_1+58632920 \\[0.5ex]
%& \qquad \quad \
%+ a_2^2\, [6075405 a_1^3+10763001 a_1^2-184769148 a_1-175937284\\[0.5ex]
%& \qquad \quad \
%-10 a_2^2 (2 (1482165 a_1-1645331)-95 a_2^2 (1311 a_1+8543)) ] \} \ne 0.
%\end{array}
$ at $(a_1,a_2)=(a_{13},a_{23})$.
This, together with the above results,
suggests that we may have parameter values such that
$ v_{3i}=0, \, i=0,1,2, \dots, 7$, $v_{38} \ne 0$,
 and so the system could have $8$ small-amplitude limit cycles,
by properly applying perturbations on the coefficients,
$ b_{013}$, $ b_{023}$, $ b_{123}$, $ b_{113}$, $ b_{302}$, $ a_{211}$, $a_1 $
and $a_2$.

Now, we want all $\varepsilon^3$-order focus values to vanish
(i.e., $M_3(h) \equiv 0$). This can be achieved by solving the coefficient
$a_{121}$ from a polynomial equation. Having obtained the conditions for
which all the $\varepsilon^1$-,  $\varepsilon^2$- and  $\varepsilon^3$-order
focus values vanish, we now use the $\varepsilon^4$-order focus values to
linearly solve for $b_{024}$, $ b_{124}$, $ b_{114}$, $ b_{303}$, $ a_{212}$
and $a_{122}$ one by one from the equations
$v_{41}=v_{42}=v_{43}=v_{44}=v_{45} = v_{46} =0$.
Then, the higher-order focus values are given by
$$
v_{47} = \tss\frac{13}{1179648}\, F_{40} F_{41}, \ \
v_{48} = \tss\frac{-13}{127401984}\, F_{40} F_{42}, \ \
v_{49} = \tss\frac{13}{244611809280}\, F_{40} F_{43},
$$
where $F_{40}$ is a common factor, and $F_{41}$, $F_{42}$ and
$F_{43}$ are polynomials in $a_1$ and $a_2$, given in Appendix A.
Similarly, we obtain the solutions of $a_1$ and $a_2$ for
$F_{41} = F_{42} = 0$, but $F_{43} \ne 0$, given as follows:
\begin{equation}
\begin{array}{ll}
%(a_1,\, a_2) = (a_{14}^i,\, \pm a_{24}^i), \quad
a_1 = a_{14}^i,\, \, a_2 = \pm a_{24}^i =\pm a_2(a_{14}^i),\quad
i=1,2, \dots, 6, \quad {\rm where} \\[0.5ex]
a_{14}^1 = -4.58252393 \cdots, \ \
a_{14}^2 = -1.72294798 \cdots, \ \
a_{14}^3 = -0.21827689 \cdots, \\[0.5ex]
a_{14}^4 = -0.09420293  \cdots, \ \
a_{14}^5 = \hspace{0.1in} 0.14811742 \cdots, \ \
a_{14}^6 = \hspace{0.1in} 1.45012903 \cdots,
\end{array}
\label{Eq_soln4}
\end{equation}
and $a_2(.)$ denotes a rational function of the variable,
which satisfy $F_{43} \ne 0$ and
$
%\det \! \left[ \begin{array}{cc}
%\frac{\partial F_{41}}{\partial a_1} &
%\frac{\partial F_{41}}{\partial a_2} \\
%\frac{\partial F_{42}}{\partial a_1} &
%\frac{\partial F_{42}}{\partial a_2} \end{array}
%\right]_{F_{11}= \tilde{F}_{12} = 0}
\det \left[\partial (F_{41},F_{42})/\partial (a_1,a_2)\right]_{F_{41}
= \tilde{F}_{42} = 0} \ne 0.
$
This suggests that with the $\varepsilon^4$-order focus values, we can obtain
$9$ small-amplitude limit cycles by properly perturbing the coefficients,
$b_{014}$, $ b_{024}$, $ b_{124}$$, $ $b_{114}$, $ b_{303}$, $ a_{212}$,
$ a_{122}$, $ a_1$ and $a_2$.

Finally, in order to have all the $\varepsilon^4$-order focus values to
become zero, we let $b_{021} = -\frac{2 a_2}{a_1^2} $. Then, we obtain
the simplified conditions, listed in Appendix B,
under which all the $\varepsilon^1$-,
$\varepsilon^2$-, $\varepsilon^3$-, $\varepsilon^4$-order focus values
vanish. Then, we use the $\varepsilon^5$-order focus values
to find $10$ small-amplitude limit cycles. Linearly solving the seven
polynomial equations, $v_{51}=v_{52} = \cdots = v_{57} = 0$ one by one for
the seven coefficients, $b_{025}$, $b_{125}$, $b_{115}$, $b_{304}$,
$a_{213}$, $a_{123}$ and $b_{022}$. Then, $v_{58}$, $v_{59}$ and $v_{510}$
are given in terms of $a_1$ and $a_2$:
$$
v_{58} = \tss\frac{187}{6193152000} \, F_{50} F_{51}, \ \
v_{59} = \tss\frac{-187}{990904320000}\, F_{50} F_{52}, \ \
v_{510} = \tss\frac{17}{11890851840000}\, F_{50} F_{53},
$$
where the common factor $F_{50}$ is a rational function of $a_1$ and $a_2$,
and $F_{5i}, \, i=1,2,3$ are polynomials in $a_1$ and $a_2$, with degrees
$6$, $7$ and $8$ with respect to $a_2^2$, respectively. $F_{51}$ and
$F_{52}$ are given in Appendix A ($F_{53}$ is omitted for brevity).
It can be shown that there are in a total $12$ real solutions for
$(a_1,a_2)$ such that $F_{51}=F_{52}=0$, but $F_{53} \ne 0$, given as follows:
\begin{equation}
\begin{split}
&a_1 = a_{15}^i,\, \, a_2 = \pm a_{25}^i =\pm a_2(a_{15}^i) \quad
i=1,2, \dots, 6, \quad {\rm where} \\
&a_{15}^1 = -2.39560267 \cdots, \ \
a_{15}^2 = -1.53681619 \cdots, \ \
a_{15}^3 = -0.38249860 \cdots,\\
&a_{15}^4 = -0.19575710 \cdots, \ \
a_{15}^5 = \hspace{0.1in} 0.05960015 \cdots, \ \
a_{15}^6 = \hspace{0.1in} 0.29402249 \cdots,
\end{split}
\label{Eq_soln5}
\end{equation}
and $a_2(.)$ denotes a rational function of the variable,
which satisfy $F_{53} \ne 0$ and
$
%\det \! \left[ \begin{array}{cc}
%\displaystyle\frac{\partial v_{58}}{\partial a_1} & \displaystyle\frac{\partial v_{58}}{\partial a_2} \vspace{0.2cm}\\
%\displaystyle\frac{\partial v_{59}}{\partial a_1} & \displaystyle\frac{\partial v_{59}}{\partial a_2}
%\end{array} \right]_{F_{51}=F_{52}=0}
\det \left[\partial (F_{51},F_{52})/\partial (a_1,a_2)\right]_{F_{51}
=F_{52}=0} \ne 0,
$
implying that we can apply perturbations on the $10$ parameters,
$b_{015}$, $b_{025}$, $b_{125}$, $b_{115}$, $b_{304}$, $a_{213}$, $a_{123}$,
$b_{022}$, $a_1$ and $a_2$ to obtain $10$ small-amplitude limit cycles
around the origin.

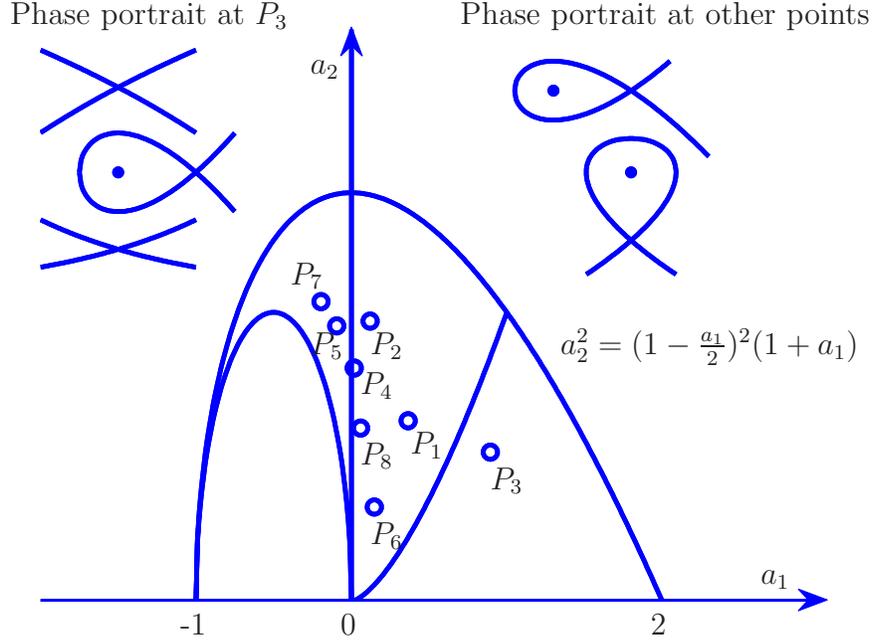
\begin{figure}[!t]
\vspace{-0.4in}
\begin{overpic}
[scale=0.9]{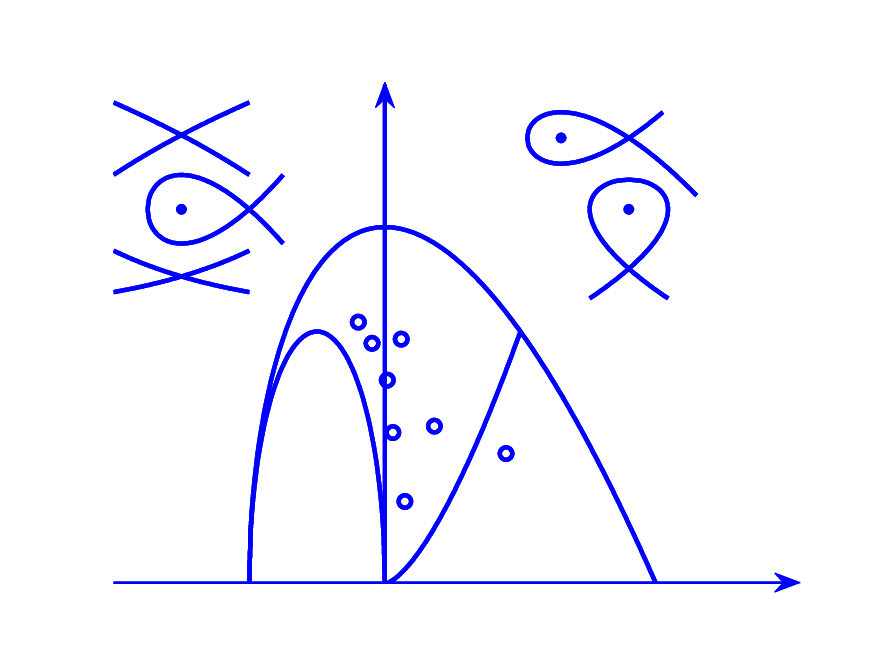}
    \put(50,23){$P_1$} \put(46,33){$P_2$} \put(58,19.5){$P_3$} \put(45,29){$P_4$}
    \put(40,33){$P_5$} \put(46,14){$P_6$} \put(38,40){$P_7$} \put(45,22){$P_8$}
    \put(27,5){-1} \put(43,5){0} \put(74,5){2} \put(85,10){$a_1$} \put(40,61){$a_2$}
    \put(10,66){Phase portrait at $P_3$} \put(55,66){Phase portrait at other points}
    \put(65,33){$a_2^2 = (1- \frac{a_1}{2})^2 (1+a_1)$}
\end{overpic}
\vspace{-0.2in}
\caption{Distribution of points $P_i$ and their corresponding phase portraits.}
\label{C2Fig1}
\end{figure}

Finally, we need to check the critical values given in equations
(\ref{Eq_soln1}), (\ref{Eq_soln2}), (\ref{Eq_soln3}),
(\ref{Eq_soln4}) and (\ref{Eq_soln5}) which are properly distributed in
the bifurcation diagram in terms of the parameters $a_1$ and
$a_2$ with the Hamiltonian function $H(x,y)$ given in (\ref{Eq7}).
(See Figure \ref{C2Fig1} in~\cite{HorozovIliev1994} for the Hamiltonian function
$H(x,y) = \frac{1}{2} (x^2+y^2) - \frac{1}{3} x^3 + a x y^2 +
\frac{1}{3} b y^3$, in terms of the parameters $a$ and $b$.)
For convenience, we define the following points
in the $a_1$-$a_2$ plane:
\begin{equation*}
\begin{split}
k=1: \, \, P_1 =&\, (\hspace{0.13in} 0.3650705869 \dots,\
                  0.4417795388 \dots) \\
k=2: \, \, P_2= &\,   (\hspace{0.13in} 0.1214887712 \dots, \
                  0.6855794168 \dots) \\
         P_3= &\,   (\hspace{0.13in} 0.8947127237 \dots, \
                  0.3648137316 \dots) \\
k=3: \, \, P_4 =&\,  (\hspace{0.13in} 0.0187162703 \dots,\
                  0.5708409903 \dots) \\
k=4: \, \, P_5 =&\, (-0.0942029335 \dots, \
                  0.6741464973 \dots) \\
         P_6 =&\,  (\hspace{0.13in} 0.1481174260 \dots, \
                  0.2303270018 \dots) \\
k=5: \, \, P_7 =&\, (-0.1957571086 \dots, \
                  0.7336772199 \dots) \\
         P_8 =&\,  (\hspace{0.13in} 0.0596001501 \dots, \
                  0.4237619510 \dots),
\end{split}
\end{equation*}
where the number $k$ denotes the order of Melnikov functions.
Note that all of these points satisfy the conditions $-1 \le a_1 \le 2$
and $ 0 \le a_2 \le (1-a_1/2)\sqrt{1+a_1}$, that is, they are inside the
curve defined by
$$
a_2^2 = \Big(1 - \frac{a_1}{2} \Big)^2 ( 1+ a_1),
$$
as shown in Figure~\ref{C2Fig1}.
But it should be noted that there are other points outside
the curve (not shown in this figure) which are also solutions.
For each $k$, there exist proper Hamiltonian functions for
which the conclusion in Theorem \ref{TheoB1} holds.
It has been seen from our solution procedures that $a_2 = 0$ is
not allowed, and none of the above cases is degenerate.
In particular, the degenerate case, defined by $a_1^3= 2 a_2^2$,
does not belong to our parameter values.
The corresponding phase portraits for the
eight sets of parameter values ($8$ points $P_i$)
are also sketched in Figure~\ref{C2Fig1}.

The above results indeed show that by using the $k$th-order Melnikov function
$M_k$, we can obtain $\lfloor\frac{4k}{3}\rfloor \!+\!4$ number small-amplitude
limit cycles bifurcating from the origin of system (\ref{Eq6}).

The proof for Theorem 1 is complete.
\end{proof}

\section{Conclusion}
\label{concl}
In this paper, we have shown that
the basis chosen in the proof of~\cite{Zoladek1995}
are not independent, leading to the conclusion of
the existence of $11$ limit cycles in this example being not true.
Further, with an example, we have demonstrated a good method combining
both advanteges of the Melnikov function method and the focus value computation
method in studying bifurcation of limit cycls. In particular,
we perturb a quadratic Hamiltonian system with cubic polynomials to
obtain $10$ small-amplitude limit cycles by using up to 5th-order Melnikov
funcstions. This illustrates the usefulness of the combination method,
and it is expected that this method can be applied to investigate other
polynomial systems to obtain more limit cycles.

\section*{Acknowledgment}
\label{acknow}
The authors thank the anonymous referees for their valuable comments
and suggestions, which have greatly helped improving the presentation
of this paper. This work was supported by the National Science and
Engineering Research Council of Canada (NSERC).

\section*{Appendix A}

\begin{equation*}
\begin{split}
F_{41} = & 37179 a_1^8-524880 a_1^7-4747248 a_1^6-12436416 a_1^5 +7737120 a_1^4
\\
& +13042944 a_1^3 -17299200 a_1^2+6945792 a_1-578816 \\
& -16 a_2^2 \big\{ 12393 a_1^6 -802548 a_1^4-102708 a_1^5
-1317600 a_1^3 -40464 a_1^2 \\
& +232128 a_1+144704
-2 a_2^2 [ 3 (11475 a_1^4-35496 a_1^3-271896 a_1^2 \\
& -38688 a_1+129712)+18088 a_2^2 (3 a_1^2+12 a_1-4-a_2^2) ]
\big\} \\
%\end{split}
%\end{equation*}
%\begin{equation*}
%\begin{split}
F_{42} = & \ 2676888 a_1^{10} \!-\! 52205877 a_1^9 \!-\!
223716978 a_1^8 \!+\! 3206239200 a_1^6 \!+\! 200795760 a_1^7\\
& -5054946912 a_1^5-3905952192 a_1^4+10386531072 a_1^3 -7205736960 a_1^2\\
&+2022961920 a_1-144704000 +2 a_2^2 \big\{ 780759 a_1^8-9325368 a_1^7 \\
&-641496672 a_1^6 -2909890656 a_1^5+558977760 a_1^4 +5294374272 a_1^3 \\
& -2824768512 a_1^2 +1366334976 a_1-486784256-16 a_2^2 \big[ 945999 a_1^6 \\
& -2537325 a_1^5-136313118 a_1^4 -397028520 a_1^3+244645056 a_1^2 \\
& +201986928 a_1-100792544 -2 a_2^2 \big( 4716225 a_1^4-23135436 a_1^3 \\
& -141272064 a_1^2+42697968 a_1+50396272+9044 a_2^2 (3495 a_1+810 a_1^2 \\
& -1682 -250 a_2^2) \big) \big] \big\} \\
\end{split}
\end{equation*}
\begin{equation*}
\begin{split}
F_{43} = & \ 5 (1366216713 a_1^{12} \!-\! 33939081972 a_1^{11} \!+\!
44272893168 a_1^{10} \!+\! 493387025040 a_1^9\\
& -628994298672 a_1^8-14032675198080 a_1^7+18889326323712 a_1^6\\
& +18007656030720 a_1^5-42717415378176 a_1^4 +28085352201216 a_1^3\\
& -6517758455808 a_1^2-261522960384 a_1+82636402688) \\
& -16 a_2^2 \big\{ 3 (3255393240 a_1^{10} -44439681807 a_1^9-361043394498 a_1^8 \\
& -1113177716208 a_1^7-333653885856 a_1^6+4370955883488 a_1^5 \\
& +418550391360 a_1^4-5262495843072 a_1^3+6238547740160 a_1^2\\
& -2585249949440 a_1+164748398080)-a_2^2 \big[ 3 (17368810155 a_1^8 \\
& -138413665080 a_1^7 -974515821120 a_1^6 -2142258103008 a_1^5 \\
& -1380949222176 a_1^4+8851316920704 a_1^3+2793260427776 a_1^2 \\
& -2239742773760 a_1-792175055104)-16 a_2^2 \big(2 (7413637185 a_1^6\\
& -9049012605 a_1^5 -244787495850 a_1^4-257911746696 a_1^3 \\
& +922435001664 a_1^2+209975885040 a_1-361224302752) \\
& +a_2^2 (3 (9653815755 a_1^4+43625458140 a_1^3-8316724720 a_1^2 \\
& -161578121840 a_1 +49510940944)-180880 a_2^2 (26034 a_1^2\\
&-28239 a_1-170778+8923 a_2^2)) \big) \big] \} \\
%\end{split}
%\end{equation*}
%\begin{equation*}
%\begin{split}
F_{51} = & \ 3365793 a_1^{12}+60938568 a_1^{11}-774250488 a_1^{10}+1966200480 a_1^9 \\
& +13136171760 a_1^8-8029124352 a_1^7-42401159424 a_1^6+61639418880 a_1^5 \\
& +11348709120 a_1^4-85053265920 a_1^3+68653025280 a_1^2-20425531392 a_1 \\
& +2343047168-8 a_2^2 \big\{ 3 (1620567 a_1^{10}+26340228 a_1^9-214842132 a_1^8 \\
& +216250560 a_1^7 +2573086176 a_1^6+131414400 a_1^5-4093628544 a_1^4 \\
& +1881934848 a_1^3 +1137593088 a_1^2-1275165696 a_1+718412800) \\
& -2 a_2^2 \big[ 3 (10180485 a_1^8+153299952 a_1^7-674144208 a_1^6-353045952 a_1^5\\
&+4636649952 a_1^4+880277760 a_1^3-3232210176 a_1^2 +170572800 a_1 \\
& +1300940032)+16 a_2^2 \big(7853517 a_1^6+134834868 a_1^5-120423348 a_1^4 \\
& -748001952 a_1^3+215457840 a_1^2+31982400 a_1-434094272 \\
& -133 a_2^2 (3 (14175 a_1^4-72216 a_1^3-415512 a_1^2-299616 a_1-611344) \\
& -8 a_2^2 (1215 a_1^2-74988 a_1-63300+8602 a_2^2)) \big) \big] \big\}
%\\
\end{split}
\end{equation*}
\begin{equation*}
\begin{split}
F_{52} = & \ 595745361 a_1^{14}+9456106860 a_1^{13}-180495550692 a_1^{12}+866884039776 a_1^{11}\\
&+1125517505040 a_1^{10}-7989977121984 a_1^9 +3366147119040 a_1^8 \\
& +34380042236928 a_1^7-59273145771264 a_1^6+7717979427840 a_1^5 \\
& +76097098183680 a_1^4-94586831216640 a_1^3+49990295040000 a_1^2 \\
& -12029752197120 a_1+1171523584000-4 a_2^2 \big\{ 908390133 a_1^{12} \\
& +9845436600 a_1^{11}-161757046008 a_1^{10}+687515327712 a_1^9 \\
& -956879159760 a_1^8-5927821906176 a_1^7+11861554007808 a_1^6\\
&+5082805251072 a_1^5-21066056398080 a_1^4+16917933938688 a_1^3 \\
& -4446179260416 a_1^2-2461820755968 a_1+1818858868736 \\
& -4 a_2^2 \big[ 3 (680430375 a_1^{10}+5494977468 a_1^9-93511317348 a_1^8 \\
& +342554370624 a_1^7 -417148084512 a_1^6-2601812555136 a_1^5\\
&+3432229497216 a_1^4 +1662404330496 a_1^3-2256035345664 a_1^2 \\
& -180750713856 a_1+499346197504) +4 a_2^2 \big(6850884663 a_1^8 \\
& +97499706480 a_1^7\!-\! 620311388976 a_1^6 \!+\! 219072133440 a_1^5
\!+\! 3577077122976 a_1^4 \\
& -2324658546432 a_1^3-1344157539072 a_1^2+1227803667456 a_1 \\
& -445085841664 +4 a_2^2 (4026591891 a_1^6+79808495508 a_1^5-88999950204 a_1^4
\\
& +243954873936 a_1^2-494926244640 a_1^3-33890758848 a_1-111271460416 \\
& -532 a_2^2 (3 (1871505 a_1^4 \!-\! 8465592 a_1^3 \!-\! 82349256 a_1^2
\!-\! 24101472 a_1 \!-\! 58663792) \\
& -4 a_2^2 (248955 a_1^2-22082340 a_1+2150500 a_2^2-13355108)))
\big) \big] \big\}
\end{split}
\end{equation*}

\section*{Appendix B}

\begin{equation*}
\begin{split}
&b_{121} = b_{112} =  b_{301} = 0, \quad
b_{031} = - \tss\frac{a_1^3+8 a_2^2}{10 a_1^2}, \quad b_{111} = -1 \\
&b_{032} = \tss\frac{12\,a_2}{25}\,  b_{022}
- \tss\frac{2 (5 a_1 + 31) a_1^6 -2 a_2^2
\,\left[ (8 a_1 + 13) a_1^3 + 4 (a_1 + 14) a_2^2 \right] }
{125 a_2 a_1^5 (a_1^3 - 2 a_2^2)} \\
&b_{122} = \tss\frac{3\,a_1}{5}\, b_{022}
- \tss\frac{2 (2 a_1 - 2) a_1^6 - a_2^2
\, \left[ (17 a_1 - 50) a_1^3+  4 (a_1 + 38) a_2^2 \right]}
{25 a_1^4 (a_1^3 - 2 a_2^2)} \\
&b_{113} = b_{302}, \quad a_{211} = - \tss\frac{2 a_2}{a_1}, \quad
a_{121} = - \tss\frac{a_1^3+8 a_2^2}{5 a_1^2} ,
\quad b_{114} = b_{303}, \quad b_{021} = - \tss\frac{2 a_2}{a_1^2} \\
&b_{023} = \tss\frac{1}{2 a_2}\, a_{123}
- \tss\frac{3+a_1}{10 a_2}\, b_{022}
- \tss\frac{(a_1 - 2)
\left\{2 a_1^6 \!+\!  a_2^2 \left[ (3 a_1 \!-\! 22) a_1^3
\!+\! 4 (a_1 \!+\! 14) a_2^2 \right] \right\} }
{25 a_2 a_1^5 (a_1^3 - 2 a_2^2)} \\
\end{split}
\end{equation*}
\begin{equation*}
\begin{split}
&b_{123} =-\,a_{213}+ \tss\frac{a_1}{2 a_2}\, a_{123}
         -\frac{a_1^2 \left[ 5 (a_1+3) a_1^3+4 (a_1-34) a_2^2 \right] }
{50 a_2} \, b_{022} \\
&\qquad \ \
- \tss\frac{
\left\{ 10 (a_1-2) a_1^9+a_2^2 \left[ (5 a_1^2-142 a_1+208) a_1^6
 -8 a_2^2 ((a_1^2-118 a_1+96) a_1^3 +8 a_2^2 (a_1-1) (a_1+21)) \right]
\right\}} {125 a_2 a_1^7 (a_1^3-2 a_2^2)} \\
&b_{302} = \tss\frac{2}{5}\, b_{022} + \tss\frac{4
        \left\{ 2 a_1^6+a_2^2 \left[ (3 a_1-22) a_1^3+4 (a_1+14) a_2^2 \right]
\right\} } {25 a_1^5 (a_1^3-2 a_2^2)} \\
&a_{212} = \tss\frac{2 a_1}{5}\, b_{022} - \tss\frac{2
       \left\{ (3 a_1+7) a_1^6-a_2^2 \left[ (3 a_1+20) a_1^3-24 (a_1+3) a_2^2
\right] \right\} } {25 a_1^4 (a_1^3-2 a_2^2)} \\
&a_{122} = \tss\frac{14 a_2}{25}\, b_{022}+ \tss\frac{2 a_2
         \left\{(15 a_1 - 32) a_1^6 - 4 a_2^2 \left[ (12 a_1 - 43) a_1^3
         + 6 (a_1 + 14) a_2^2 \right] \right\} }
{125 a_1^5 (a_1^3 - 2 a_2^2)} \\
&b_{024} = \tss\frac{a_1-2}{2 a_1 a_2}\, a_{213} - \tss\frac{
(3 a_1 - 1) a_1^3 - 4 (9 a_1 - 13) a_2^2 } {20 a_2^2 (a_1^3-4 a_2^2)}\, a_{123}
- \tss\frac{4 a_1^2 (a_1-2)}{25 (a_1^3-4 a_2^2)}\, b_{022}^2 \\
&\qquad \ \, + \tss\frac{1}{500 a_1^3 a_2^2 (a_1^3-2 a_2^2) (a_1^3-4 a_2^2)}
        \big\{ (15 a_1^2+40 a_1-15) a_1^9
-2 a_2^2 \big[ (113 a_1^2 \!+\! 366 a_1 \\
& \qquad \ \,
-609) a_1^6 \!-\! 4 a_2^2 \big((16 a_1^2 \!+\! 571 a_1 \!-\! 781) a_1^3
\!-\! 8 a_2^2 (3 a_1^2 \!+\! 102 a_1 \!-\! 166)\!
\big) \! \big]\! \big\} b_{022} \\
& \qquad \ \,
       + \tss\frac{2}{2500 a_1^8 a_2^2 (a_1^3-2 a_2^2) (a_1^3-4 a_2^2)}
       \big\{ (30 a_1^2-70 a_1+20) a_1^{12}
-a_2^2 \big[ (65 a_1^3+857 a_1^2
\\
& \qquad \ \, -864 a_1 \!-\!1420) a_1^9 \!-\! 8 a_2^2 \big((9 a_1^3 \!+\!
849 a_1^2 \!-\! 314 a_1 \!-\! 1640) a_1^6 \!-\! 2 a_2^2 ((18 a_1^3 \\
& \qquad \ \,
+1759 a_1^2 \!-\! 740 a_1 \!-\! 2500)
a_1^3 \!-\! 8 a_2^2 (11 a_1^3 \!+\! 207 a_1^2 \!-\! 108 a_1
\!+\! 100)) \big) \big]  \big\} \\
&b_{303} = -\, \tss\frac{2}{a_1}\, a_{213}
          + \tss\frac{3 (a_1^3-12 a_2^2)}{5 a_2 (a_1^3-4 a_2^2)}\,  a_{123}
        + \tss\frac{16 a_1^2 a_2}{25 (a_1^3-4 a_2^2)}\, b_{022}^2 \\
& \qquad \ \,
        - \tss\frac{1}{125 a_1^3 a_2 (a_1^3-4 a_2^2) (a_1^3-2 a_2^2)}
         \big\{ 15 (a_1+3) a_1^9-2 a_2^2 \big[ (113 a_1+417) a_1^6 \\
& \qquad \ \, -16 a_2^2 ((4 a_1+107) a_1^3 -2 (3 a_1+83) a_2^2)
\big] \big\} b_{022} \\
& \qquad \ \,
        - \tss\frac{2}{625 a_1^8 a_2 (a_1^3-4 a_2^2) (a_1^3-2 a_2^2)}
        \big\{ 30 (a_1-2) a_1^{12}  -a_2^2
\big[ (65 a_1^2 \!+\! 1062 a_1 \!-\! 40) a_1^9 \\
& \qquad \ \,
-8 a_2^2 \big((9 a_1^2+967 a_1+370) a_1^6
-4 a_2^2 ((9 a_1^2 +935 a_1+450) a_1^3 \\
& \qquad \ \,  -4 a_2^2 (11 a_1^2+204 a_1-50)) \big) \big] \big\}
\end{split}
\end{equation*}
\begin{equation*}
\begin{split}
&b_{124} = \tss\frac{a_1^3 \!\!-\! 2 a_1^2 \!\!-\! 4 a_2^2}{2 a_1^2 a_2}\,
a_{213} - \tss\frac{125
       \left\{ (3 a_1 -1) a_1^6 - 4 a_2^2 \left[ (12 a_1 - 5) a_1^3
- 4 (9 a_1 + 14) a_2^2 \right] \right\} }
{2500 a_1^2 a_2^2 (a_1^3-4 a_2^2)}\, a_{123} \\
& \qquad \ \,
- \tss\frac{2 [ (2 a_1+1) a_1^3-4 (2 a_1+7) a_2^2 ] }{25 (a_1^3-4 a_2^2)}
\, b_{022}^2 + \tss\frac{1}{500 a_1^5 a_2^2 (a_1^3-2 a_2^2) (a_1^3-4 a_2^2)} \\
& \qquad \ \, \times
       \big\{ 5 (3 a_1^2 + 8 a_1 - 3) a_1^{12}
         - 2 a_2^2 \big[ (53 a_1^2 + 436 a_1 - 459) a_1^9 \\
&\qquad \ \, +4 a_2^2 ((26 a_1^2-1141 a_1+624) a_1^6
         -4 a_2^2 ((53 a_1^2-758 a_1+474) a_1^3 \\
&\qquad \ \, -16 a_2^2 (a_1^2 - 48 a_1 + 21))) \big] \big\} b_{022}
       +  \tss\frac{1}{1250 a_1^9 a_2^2 (a_1^3 - 2 a_2^2) (a_1^3 -
        4 a_2^2)} \\
& \qquad \ \ \, \times
       \big\{ 10 (3 a_1^2 - 7 a_1 + 2) a_1^{14}
-a_2^2 \big[ (65 a_1^3+617 a_1^2-904 a_1-2140) a_1^{11} \\
& \qquad \ \,
          -4 a_2^2 \big( (218 a_1^3 + 1205 a_1^2 + 994 a_1 -5000) a_1^8
          -4 a_2^2 ((156 a_1^3 \\
& \qquad \ \,
+1731 a_1^2+934 a_1-4140) a_1^5
-4 a_2^2 ((15 a_1^3+1392 a_1^2
+54 a_1+20) a_1^2 \\
& \qquad \ \,
-48 a_2^2 (a_1^2+20 a_1-56))) \big) \big]
\big\}
\end{split}
\end{equation*}


%% References


%\bibliographystyle{elsarticle-num}
%\bibliography{TY_JDE_bib}

\end{document}